\documentclass[11pt]{article}
\usepackage{amsmath}
\usepackage{amssymb}
\usepackage{amsfonts}
\usepackage{graphicx}
\usepackage{epsfig}
\usepackage{color}
\usepackage{epsfig,amsmath,amssymb,latexsym,color,placeins}
\usepackage{algorithmic}
\usepackage[ruled,vlined]{algorithm2e}
\usepackage{graphicx}
\usepackage{subcaption}
\usepackage{xcolor}
\usepackage{epstopdf}
\usepackage[toc,page]{appendix}
\usepackage{siunitx,array,multirow}

\def\bfc{{\bf c}}

\newtheorem{theorem}{Theorem}
\newtheorem{lemma}{Lemma}

\newtheorem{example}{Example}

\newenvironment{proof}{\begin{trivlist}\item[]{\emph{Proof.}}}
               {\hfill$\Box$\end{trivlist}}

\begin{document}
\title{ A Spline-Based Collocation Method for \\Stokes and Navier-Stokes equations}
\author{Jinsil Lee
\footnote{jl74942@ewha.ac.kr, Department of Mathematics, Ewha W. University, Seoul 120-750, South Korea,This work was funded by NRF-2019R1A6A1A11051177.}}
\date{}
\maketitle

\begin{abstract}
	In the paper, we propose a collocation method based on multivariate polynomial splines over triangulation or tetrahedralization for solving Stokes and Navier-Stokes equations.  
	We start with a detailed explanation of the method for the Stokes equation and then extend the study to the Navier-Stokes equations. And we provide the analysis for our method and show the existence of the spline approximation. We shall show that the numerical solution can approximate the exact PDE solution very well over several domains. Then we present several numerical experimental results to demonstrate the performance of the method over the 2D and 3D settings. Also, we apply the IPBM method to our method to find the solution over several curved domains effectively. In addition, we present a comparison with the existing multivariate spline methods in \cite{AL02} and several existing methods to show that the new method produces a similar and sometimes more accurate approximation in a more efficient fashion.  
\end{abstract}



\section{Introduction}
The Navier-Stokes equation is widely used to describe various physical phenomena, such as wave breaking in the ocean and airflow around airplane wings. In steady flows, viscous transport dominates over convection, which can be neglected when the fluid moves slowly enough. This simplifies the momentum equation of the general Navier-Stokes equations into a linear equation, resulting in the system \eqref{stokes}. The Stokes equations have garnered substantial attention from researchers due to their close relationship with the Navier-Stokes equations. Especially, the Stokes and Navier-Stokes equations with variable viscosity are crucial for modeling creeping flows observed in various geophysical contexts, as demonstrated by prior studies such as \cite{H83}, \cite{M84}, \cite{MZG96}, and \cite{ZYM07}. These equations, central to predicting phenomena as varied as glacial motion in the face of climate change and the aerodynamic design of transportation, demand solutions of the highest precision and adaptability. In \cite{JKN15}, researchers applied the finite element method to analyze incompressible Stokes equations with a spatially varying viscosity parameter $\mu(x)\geq \mu_{\min}>0$.

In this paper, we propose and investigate a novel collocation method based on multivariate splines for numerically solving the Stokes equations over the domain in $\mathbb{R}^d$ with $d\ge 2$. Let us consider the Navier-Stokes equation:
\begin{align}
	\label{nstokes}
	-\mu\Delta \textbf{u}+\textbf{u}\cdot \nabla \textbf{u}+\nabla p &=f ~~~\text{in}~\Omega \subseteq \mathbb{R}^d \\
	\nabla \cdot \textbf{u} &=0 ~~~\text{in}~\Omega \\
	\textbf{u}&=g ~~~\text{on}~\partial \Omega
\end{align}
where the unknowns velocity $\textbf{u}=(u_1,\cdots, u_d)^T$ of the fluid and pressure $p$ and $\mu$ is the kinematic viscosity, $f=(f_1, \cdots,f_d)$ represents the externally applied forces and $g=(g_1, \cdots, g_d)$ the velocity at the boundary. 
Numerical solutions to these Stokes equations have been recently studied extensively. There are many different numerical methods that have been developed such as conforming/nonconforming finite element methods, weak Galerkin method,  finite volume methods, discontinuous Galerkin finite element methods, virtual element methods, WEB-spline method(Weighted extended B-spline method), etc.

The finite element method is a popular approach for solving the Stokes equations, as documented in \cite{CR73}, \cite{GR86}, and \cite{MY17}. However, one of the primary challenges in solving the Stokes and Navier-Stokes equations is the coupling of velocity and pressure variables in a saddle-point system. The velocity error is pressure-dependent, which becomes particularly significant when dealing with small viscosity $\mu$. This issue can lead to substantial errors or singularities in pressure that affect velocity simulation due to the coupling between velocity and pressure. To address this challenge, numerous methods have been developed, including divergence-free finite element methods discussed in \cite{FN13}, \cite{SV85}, \cite{Z11}, \cite{M20}, \cite{MWY15}, and \cite{MYZ21}. These methods explicitly construct a divergence-free basis and approximate velocity from weakly or exactly divergence-free subspaces. In \cite{MYZ21}, enhanced discontinuous Galerkin (DG) finite element methods are used to solve Stokes equations in the primary velocity-pressure formulation, achieving pressure-robustness and establishing optimal-order error estimates.
Researchers in \cite{BMW19} utilized the weak Galerkin finite element method for the Stokes equations. Additionally, in \cite{KKD06}, a mesh-free method that employs weighted extended B-splines as shape functions was employed to solve the Stokes problem, accompanied by convergence and condition number estimates.

In this paper, our approach focuses on utilizing the spline collocation method, augmented with multivariate splines, to address both Stokes and Navier-Stokes equations. Beyond the primary advantages offered by spline functions, such as their degree adaptability, customizable smoothness, and the intrinsic property of Bernstein-Bézier polynomials' unity partition, our approach offers several additional benefits compared to traditional methods like finite element methods, discontinuous Galerkin methods, and virtual element methods.

For instance, our method eliminates the need for a weak formulation of the partial differential equation (PDE) solution, thereby obviating the requirement for numerical quadrature during computation. Additionally, we possess greater flexibility in managing discontinuities stemming from PDE coefficients by readily adjusting the positions of collocation points on either side of discontinuous curves or surfaces. Furthermore, when dealing with PDEs over curved domains, we can easily reposition boundary points to determine solutions instead of resorting to increased mesh density near the boundary or integral computations along it.

Moreover, our multivariate spline-based collocation approach offers the capacity to enhance approximation accuracy by increasing the degree, a more cost-effective alternative to pursuing a uniform refinement of the underlying triangulation or tetrahedralization within the confines of a computer's memory budget.
We will present a series of numerical tests in 2D and 3D to highlight the effectiveness of our method. Our objectives are to evaluate its performance under various conditions, including:
\begin{itemize}
	\item[1.] 2D Stokes equation with a continuous and discontinuous viscosity $\mu$ 
	\item[2.] 3D Stokes equation with a continuous and discontinuous viscosity $\mu$ 
	\item[3.] 2D Navier-Stokes equation with a continuous and discontinuous viscosity $\mu$ 
	\item[4.] 2D Stokes equation and Navier-Stokes equation over the curved domain.
\end{itemize}
Our goal is to demonstrate that the spline-based collocation method is versatile and adept, capable of tackling the inherent complexities of these mathematical challenges.

The existence and uniqueness results for the weak formulation of \eqref{stokes} with Dirichlet boundary conditions are given in \cite{GR86} as follows:
\begin{theorem}
	Let $\Omega$ be a bounded domain in $\mathbb{R}^d$ with Lipschitz boundary $\partial \Omega$ and let $f\in H^{-1}(\Omega)$ and $g\in H^{\frac{1}{2}}(\partial \Omega)$ satisfying
	$$\int_{\partial \Omega} g\cdot n=0,$$ 
	where $n$ is a unit normal vector.
	There is a unique solution $(u,p)\in H^1(\Omega)^d\times L^2_0(\Omega)$ of 
	\begin{align*}
		-\mu\Delta \textbf{u}+\nabla p &=f ~~~\text{in}~\Omega \subseteq \mathbb{R}^d \\
		\nabla \cdot \textbf{u} &=0 ~~~\text{in}~\Omega \\
		\textbf{u}&=g ~~~\text{on}~\partial \Omega
	\end{align*}
\end{theorem}
It is possible to use $H^2(\Omega)$ as we can use spline
functions in $H^2(\Omega)$ if the degree of splines is sufficiently large.  Indeed, how to use such spline functions has been explained in \cite{ALW06},\cite{LL22},\cite{LL23},\cite{L23} and \cite{S15}. 
Since the equations involve $\nabla p$, the pressure is determined up to an additive constant. To have uniqueness, one can require it to have zero mean. We first introduce $\hat{H}^k(\Omega)=\{ p\in H^k(\Omega): \int_\Omega p =0\}$. And we define a norm $\|(\textbf{u}, p)\|_{H^k}$ on $\tilde{H}^k_0=[H^k(\Omega)\cap H^1_0(\Omega)]^d\times  \hat{H}^{k-1}(\Omega)$ which is defined by 
$$\|(\textbf{u}, p)\|_{H^k}=\sum_{i=1}^d \| u_i\|_{H^k(\Omega)}+\| p\|_{H^{k-1}(\Omega)}.  $$
And we denote $\tilde{H}_0=\tilde{H}^2_0.$
To establish the convergence of the collocation solution $s$ as the size of $\triangle$ goes to zero, we define a new norm $\|(u_1,u_2, p)\|_L$ on $\tilde{H}_0$ as follows. 
\begin{equation}
	\label{stokesH2norm2}
	\| (\textbf{u}, p ) \|_{L}=\| \mu \Delta \textbf{u}+\nabla p\|_{L^2(\Omega)}+\| \nabla \cdot \textbf{u} \|_{L^2(\Omega)}
\end{equation}
We mainly show that the new norm is equivalent to the standard norm on $\tilde{H}_0$.   
\begin{theorem}
	\label{mjlai031820211stokes} Let $\mu$ be a positive constant.
	There exist two positive constants $A$ and $B$ such that 
	\begin{equation}
		\label{equivalentstokes}
		A \|\textbf{u} \|_{H} \le \|\textbf{u} \|_L \le B \| \textbf{u}\|_{H}, \quad \forall \textbf{u}\in \tilde{H}_0.
	\end{equation}
\end{theorem}
See the proof of Theorem~\ref{mjlai031820211stokes} in a later section. Letting $(u_1, u_2, p)\in \tilde{H}_0$
be the solution of (\ref{Stokes2}) and $(u_{s_1}, u_{s_2}, u_{s_p})$ be the spline solution of (\ref{Stokes2}), we
use the above inequality to have
$$
A\|(u_1-u_{s_1}, u_2-u_{s_2},p- u_{s_p})\|_{H}\le \|(u_1-u_{s_1}, u_2-u_{s_2},p- u_{s_p})\|_L.
$$ 
When $\|(u_1-u_{s_1}, u_2-u_{s_2},p- u_{s_p})\|_L $ is small, our Theorem~\ref{mjlai031820211stokes}
implies that $\|(u_1-u_{s_1}, u_2-u_{s_2},p- u_{s_p})\|_{H}=\| u_1-u_{s_1}\|_{H^2(\Omega)} +\| u_2-u_{s_2}\|_{H^2(\Omega)}+\| p-u_{s_p}\|_{H^1(\Omega)}$ is small.
Our paper is structured as follows: Section 2 covers the introduction of useful properties about the weak formulation of the Navier-Stokes equation and multivariate splines and their attributes. In Section 3, we delineate the spline collocation methods applied to Stokes and Navier-Stokes equations and assert the convergence of our approach. Section 4 presents numerical outcomes from multiple 2D and 3D examples, showcasing our method's efficacy. The final section unveils the IPBM method, displaying its performance on various curved domains. Throughout, we compare our technique’s precision against existing numerical methods to underscore its accuracy.

\section{Preliminaries}
The weak formulation of the problem \eqref{nstokes} is given as follows: 
\begin{align}
	\label{weaknstokes}
	\mu \int_\Omega \nabla  \textbf{u}\nabla  \textbf{v}+b(\textbf{u},\textbf{u},\textbf{v})+\int_\Omega (\nabla \cdot \textbf{v}) p &=\int_\Omega f\textbf{v} ~~~\forall \textbf{v}\in [H^1_0(\Omega)]^d,\\
	\int_\Omega (\nabla \cdot \textbf{u}) q &=0 ~~~\forall q\in \hat{H}^0(\Omega),
\end{align}
where the trilinear form $b$ is defined as
\begin{align*}
	b(\textbf{u},\textbf{v},\textbf{w})=\int_\Omega (\textbf{u}\cdot \nabla) \textbf{v}\textbf{w}+\frac{1}{2}(\nabla \cdot \textbf{u}) \textbf{v}\textbf{w}.
\end{align*}
Then we have the following properties \cite{GR86},\cite{T84}:
$$b(\textbf{u},\textbf{v},\textbf{w})=\int_\Omega \frac{1}{2}(\textbf{u}\cdot \nabla) \textbf{v}\textbf{w}-\frac{1}{2}(\nabla \cdot \textbf{u}) \textbf{w}\textbf{v} ~~\forall \textbf{u},\textbf{v},\textbf{w}\in [H^1_0(\Omega)]^d,$$
and
\begin{align}\label{operbineq}
	|b(\textbf{u},\textbf{v},\textbf{w})|\le C_1\|\nabla \textbf{u}\|\|\nabla \textbf{v}\|\|\nabla \textbf{w}\|~~\forall \textbf{u},\textbf{v},\textbf{w}\in [H^1_0(\Omega)]^d
\end{align}
where
$$C_1=\sup_{\textbf{u},\textbf{v},\textbf{w}\in [H^1_0(\Omega)]^d}\dfrac{|b(\textbf{u},\textbf{v},\textbf{w})|}{\|\nabla \textbf{u}\|\|\nabla \textbf{v}\|\|\nabla \textbf{w}\|}$$
is a positive constant depending only on the domain $\Omega.$

Next, we quickly introduce bivariate spline functions. Before we start, we first review some facts about triangles. Given a triangle $T$, we write $|T|$ for the length of its longest edge and $\rho_T$ for the radius of the largest disk that can be inscribed in $T$. 
For any polygonal domain $\Omega\subset \mathbb{R}^d$ with $d= 2$, let $\triangle:=\{T_1,\cdots, T_n\}$ be a triangulation of 
$\Omega$ which is a collection of triangles and $\mathcal{V}$ be the set of vertices of $\triangle$. {We called a triangulation as a quasi-uniform triangulation if all triangles $T$ in $\triangle$ have comparable sizes in the sense that $\frac{|T|}{\rho_T}\le C<\infty,~\text{for all triangles} ~ T\in \triangle,$
	where $\rho_T$ is the inradius of $T$. Let $|\triangle|$ be the length of the longest edge in $\triangle.$} For a triangle 
$T=(v_1,v_2,v_3) \in \Omega,$ we define the barycentric coordinates $(b_1, b_2,b_3)$ 
of a point $(x,y)\in \Omega$. These 
coordinates are the solution to the following system of equations 
\begin{eqnarray*}
	b_1+b_2+b_3=1\\
	b_1 v_{1}+b_2 v_{2}+b_3 v_{3}=(x,y)
\end{eqnarray*}
and are nonnegative if $(x,y)\in T.$ We use the barycentric coordinates to define the Bernstein polynomials of degree $D$:
\begin{eqnarray*}
	B^T_{i,j,k}(x,y):=\frac{{D}!}{i!j!k!}b_1^i b_2^j b_3^k, ~i+j+k=D,
\end{eqnarray*}
which form a basis for the space $\mathcal{P}_D$ of polynomials of degree $D$. 
Therefore, we can represent all $s\in \mathcal{P}_D$ in B-form:
\begin{eqnarray*}
	s|_T=\sum_{i+j+k=D}c_{ijk}B^T_{ijk}, \forall T\in \triangle,
\end{eqnarray*}
where the B-coefficients $c_{i,j,k}$ are uniquely determined by $s$. 

Moreover, for given $T=(v_1,v_2,v_3)\in \triangle$, 
we define the associated set of domain points to be 
\begin{equation}
	\label{domainpoints}
	\mathcal{D}_{D,T}:= \{\frac{iv_{1}+jv_{2}+kv_{3}}{D} \}_{i+j+k=D}.
\end{equation}
Let $\mathcal{D}_{D,\triangle} = \cup_{T\in \triangle} \mathcal{D}_{D,T}$ be the domain 
points of triangulation $\triangle$ and degree $D$.  

We use the discontinuous spline space $S^{-1}_D(\triangle):= 
\{s|_{T} \in \mathcal{P}_D, T\in \triangle\}$ as a base. 
Then we add the smoothness conditions to define the space 
$\mathcal{S}^r_D:=C^r(\Omega)\cap S^{-1}_D(\triangle).$ 
Let {$\bfc$} be the coefficient vector of $s\in S^{-1}_D(\triangle)$ and $H$ be the matrix that consists of 
the smoothness conditions across each interior edge of $\triangle$.  
It is known that $H{\bfc}=0$ if and only if $s\in C^r(\Omega)$ 
(cf. \cite{LS07}).   

Similarly, for trivariate splines, let $\Omega\subset\mathbb{R}^3$ and $\triangle$ be a tetrahedralization of $\Omega$. Trivariate splines can be defined in a manner analogous to the 2D splines previously described.
A lot of literatures including \cite{ALW06},  \cite{LL22}, \cite{LL23}, \cite{LW04}, \cite{L23},\cite{S15}, and \cite{S19}, offers comprehensive insights into the utilization of these spline functions.

We suppose that $Q$ is the quasi-interpolation operator. Then, it has the following properties
\begin{lemma}([Lai and Schumaker, 2007\cite{LS07}])
	\label{lem1}
	There exists a constant $K$ such that if $u \in H^{m+1}(\Omega)$ for some $0\le m\le d$, then
	$$\|D^\alpha (u-Qu)\|_{L^2(\Omega)}\le K |\triangle|^{m+1-|\alpha|}|u|_{m+1,2,\Omega},$$
	for all $0\le |\alpha|\le m.$ If $\Omega$ is convex, the constant $K$ depends only on $d, \triangle$. If $\Omega$ is not convex, $K$ also depends on the Lipschitz constant of the boundary of $\Omega$.
\end{lemma}
As mentioned in \cite{A15}, by Markov inequality in \cite{LS07}, we can get the inverse estimates
\begin{lemma}
	\label{lem2}
	There exists a constant $K$ such that if $0\le l\le s$ and $1\le p,q \le \infty$, then
	$$\|s\|_{s,p}\le K |\triangle|^{l-s+\min(0,\frac{d}{p}-\frac{d}{q})}\|s\|_{l,q},$$
	for all $s\in \mathcal{S}^r_D$.
\end{lemma}
Moreover, from Lemma \ref{lem1}, we can get the following inequality
\begin{align}\label{prelimeq3}
	\|Qu\|_{k,p}\le K\|u\|_{k,p}
\end{align}
for $1\le p\le \infty$ and $0\le k\le D.$

Also, we define a quasi-interpolation spline $Qu$ 
satisfying
$$\int_T (\nabla \cdot u) v =\int_T (\nabla \cdot Qu) v, ~~\forall v\in [P_D(\Omega)]^d$$
and
$$\|D^\alpha (u-Qu)\|_{L^2(\Omega)}\le K |\triangle|^{m+1-|\alpha|}|u|_{m+1,2,\Omega},$$
for all $0\le |\alpha|\le m.$ If $\Omega$ is convex, the constant $K$ depends only on $d, \triangle$.

\section{Our Proposed Algorithms and Their Convergence Analysis}
\subsection{Spline Collocation Methods for Stoke and Navier-Stokes Equations}
For the sake of clarity, we will utilize the given equation to elucidate our method. Our objective is to demonstrate that the numerical solution converges to the true solution. For convenience, we consider the 2D case. Let $\mu>0$ be a positive function. As previously mentioned, we will employ $C^1$ spline functions. The domain points, as referenced in \cite{LS07} and further detailed in the subsequent section, will serve as our collocation points.

For brevity, let $s$ represent a $C^1$ spline of degree $D$ defined over a triangulation $\triangle$ of $\Omega$. Let $\xi_i, i=1, \cdots, N$ be the domain points of $\triangle$ and $D'$, noting that $D'$ might differ from $D$. Our spline-based collocation approach aims to find a spline function $s_1, s_2, s_p$ that satisfies the following equation 
\begin{equation}
	\label{Stokes2}
	\begin{cases}
		-\mu(\xi_i) \Delta \textbf{s}(\xi_i) +\nabla s_p(\xi_i) & = f(\xi_i), \quad \xi_i\in \Omega\subset \mathbb{R}^d, \cr 
		\nabla \cdot~ \textbf{s}(\xi_i) & = 0,  \quad \xi_i\in  \Omega \cr
		\textbf{s}(\xi_i) & = g(\xi_i),  \quad \xi_i\in \partial \Omega 
	\end{cases}
\end{equation}
where $\textbf{s}=(s_1,s_2)^T$.
A spline space is a linear vector space. We mainly look for the coefficient vector of $s$ satisfies (\ref{Stokes2}).  Clearly, (\ref{Stokes2}) leaves a linear system that may not have a unique solution. It may be an over-determined linear system if $D'\ge D$ or an underdetermined linear system if $D'<D$.  Our method is to use a least squares solution if the system is overdetermined or a sparse solution if the system is underdetermined (cf. \cite{LW21}).   
Let us first explain a collocation method based on bivariate splines/trivariate splines 
for a solution of Stokes equation:
\begin{equation}
	\label{stokes}
	\begin{cases}
		-\mu\Delta \textbf{u}+\nabla p &=f ~~~\text{in}~\Omega \subseteq \mathbb{R}^d \\
		\nabla \cdot \textbf{u} &=0 ~~~\text{in}~\Omega \\
		\textbf{u}&=g ~~~\text{on}~\partial \Omega  
	\end{cases}
\end{equation}
where $\Delta$ is a standard linear differential 
operator.   
We want to find a numerical solution of the PDE in \eqref{stokes} based on bivariate/trivariate 
spline functions which were explained in the previous section.

For given $\triangle$ be a triangulation, we choose a set of domain points 
$\{ \xi_i\}_{i=1,\cdots, N}$ explained in the previous section 
as collocation points 
and find the coefficient vector $\textbf{c}_1, \textbf{c}_2$ and $\textbf{c}_p$ of spline functions $\displaystyle s_1 =\sum_{t\in \triangle}\sum_{|\alpha|=k}c^t_{1,\alpha} \mathcal{B}^t_\alpha , \displaystyle s_2 =\sum_{t\in \triangle}\sum_{|\alpha|=k}c^t_{2,\alpha} \mathcal{B}^t_\alpha, \displaystyle s_p =\sum_{t\in \triangle}\sum_{|\alpha|=k}c^t_{p,\alpha} \mathcal{B}^t_\alpha$
satisfying the following equation at those points
\begin{equation}
	\label{PDE2stokes}
	\begin{cases} -\mu \Delta \textbf{s}(\xi_i) +\nabla s_p(\xi_i) & = f(\xi_i), \quad \xi_i\in \Omega\subset \mathbb{R}^d, \cr 
		\nabla \cdot~ \textbf{s}(\xi_i) & = 0,  \quad \xi_i\in  \Omega \cr
		\textbf{s}(\xi_i) & = g(\xi_i),  \quad \xi_i\in \partial \Omega 
	\end{cases}
\end{equation}
where $\textbf{s} =(s_1, s_2)^T$, $\{ \xi_i\}_{i=1,\cdots, N} \in \mathcal{D}_{D',\triangle}$  are the domain points of $\triangle$ of degree $D$
as explained in (\ref{domainpoints}) in the previous section. 
Using these points, we have the following matrix equation:
\begin{align*}
	\begin{bmatrix}-\mu \Delta(\mathcal{B}^t_\alpha (\xi_i))\end{bmatrix}
	\textbf{c}_1+\begin{bmatrix}(\mathcal{B}^t_\alpha (\xi_i))_x\end{bmatrix}\textbf{c}_p=[f_1(x_i, y_i)]\\
	\begin{bmatrix}- \mu\Delta(\mathcal{B}^t_\alpha (\xi_i))\end{bmatrix}
	\textbf{c}_2+\begin{bmatrix}(\mathcal{B}^t_\alpha (\xi_i))_y\end{bmatrix}\textbf{c}_p=[f_2(x_i, y_i)]\\
	\begin{bmatrix}(\mathcal{B}^t_\alpha (\xi_i))_x\end{bmatrix}\textbf{c}_1+\begin{bmatrix}(\mathcal{B}^t_\alpha (\xi_i))_y\end{bmatrix}\textbf{c}_2=0
\end{align*}
And we can rewrite it as the saddle point problem:
\begin{align*}
	K\textbf{c} := \begin{bmatrix} -\mu \Delta(\mathcal{B}^t_\alpha (\xi_i)) &O&
		(\mathcal{B}^t_\alpha (\xi_i))_x\\
		O&-\mu \Delta(\mathcal{B}^t_\alpha (\xi_i)) &(\mathcal{B}^t_\alpha (\xi_i))_y\\
		(\mathcal{B}^t_\alpha (\xi_i))_x&(\mathcal{B}^t_\alpha (\xi_i))_y &O
	\end{bmatrix}\begin{bmatrix}\textbf{c}_1\\ \textbf{c}_2\\ \textbf{c}_p\end{bmatrix}=\begin{bmatrix} f_1(x_i, y_i)\\f_2(x_i, y_i)\\O\end{bmatrix}=\textbf{f}\\
\end{align*}
where $O$ is a zero matrix. 
Furthermore, the coefficients are constrained by $H_1\textbf{c}=0$, reflecting the smoothness conditions inherent to $\mathcal{S}$. Here, the coefficient vector is represented as $\textbf{c}=(\textbf{c}_1, \textbf{c}_2, \textbf{c}_p)^T$. The boundary condition is given by $B\textbf{c}=G$, with matrices defined as $B=\begin{bmatrix}B_1& O& O\\ O& B_1&O \end{bmatrix}$ and $G=\begin{bmatrix}G_1\\G_2 \end{bmatrix}$.

The system can be represented as:
$$\begin{bmatrix}
	K\\B\\H
\end{bmatrix}\begin{bmatrix}\textbf{c}_1\\ \textbf{c}_2\\ \textbf{c}_p\end{bmatrix}=\begin{bmatrix}\textbf{f}\\ G\\ 0\end{bmatrix}.$$ This system can be solved using the backslash operator in MATLAB.

To solve the Navier-Stokes equations, we employ Newton's method in Algorithm \ref{alg:NSequation}.
\begin{algorithm}
	\caption{Algorithm for Navier-Stokes Equations}\label{alg:NSequation}
	\begin{algorithmic}[0]
		\STATE{\textbf{Initialization:} For given parameters $\mu$ and $\epsilon$, set a large value $\mu^0$ such that $\mu^0>\mu$. For our implementation, we set $\mu^0=1$.}
		
		\STATE{\textbf{Step 1}: Determine the initial value $u_0\in \mathcal{S}^r_D$ via the Stokes equation:
			\begin{align*}
				-\mu^0 \Delta u_0+\nabla p_0&=f ~~\text{in}~ \Omega \\
				\nabla \cdot u_0=0 ~~\text{in}~ \Omega &~~
				u_0=g,  ~~\text{on}~ \partial \Omega 
		\end{align*}}
		
		\STATE{\textbf{Step 2}: For given spline approximation $u_{n}\in \mathcal{S}^r_D$, compute $u_{n+1}\in \mathcal{S}^r_D$ satisfying:
			\begin{align*}
				-\mu^0 \Delta u_{n+1}+u_{n}\cdot \nabla u_{n+1}+u_{n+1}\cdot \nabla u_{n}+\nabla p_{n+1}&=f+u_{n}\cdot \nabla u_{n} ~~\text{in}~ \Omega \\
				\nabla \cdot u_{n+1}&=0 ~~\text{in}~ \Omega\\
				u_{n+1}&=g,  ~~\text{on}~ \partial \Omega 
			\end{align*}
			Use a similar approach as in solving the Stokes equation.}

		\STATE{ \textbf{Step 3}:
			If $\|u_n-u_{n+1}\|<\epsilon$, Terminate. Otherwise, execute Step 1 and Step 2 for $\mu^*$. If it converges, stop.}
		
		\STATE{\textbf{Step 4}: Select a new $\mu^1$ such that $\mu^0>\mu^1>\mu^*$ and iterate through Steps 1-3.}
	\end{algorithmic}
\end{algorithm}

\subsection{Convergence Analysis of the Stokes Equation}
Next, we consider the numerical solution $(u_{s_1}, u_{s_2}, u_{s_p})$ for equation (\ref{PDE2stokes}). 
We aim to demonstrate that as we refine the mesh, reducing the size of $|\triangle|$, the numerical solution converges to the exact solution $(u_1, u_2, p)$, i.e.,
$\|(u_1, u_2, p)-(u_{s_1}, u_{s_2}, u_{s_p})\|_{H}<\epsilon_{\triangle} $.
\begin{lemma}
	There exists a positive constant $C_\beta$ independent of $|\triangle|$ such that all $w\in W:= L^2_0(\Omega) \cap \mathcal{S}^r_D(\triangle)$ and $|\triangle|$ small enough,
	\begin{align}
		\sup_{v\in \mathcal{S}^r_D(\triangle)} \dfrac{\int_\Omega \nabla \cdot v w}{\|\nabla v\|}\geq C_\beta \|w\|.
	\end{align}
\end{lemma}
Given any $w\in W \subseteq L^2_0(\Omega)$, reference \cite{GR86} establishes the existence of a function $\hat{v}\in H^1_0(\Omega)$ satisfying
\begin{align}\label{infsup_inequpre}
	\dfrac{\int_\Omega \nabla \cdot \hat{v} w}{\|\nabla \hat{v}\|}\geq C_0 \|w\|
\end{align}
for a specific constant $C_0$. 
Let $s_{\hat{v}}$ be a quasi-interpolatory spline in $\mathcal{S}^r_D(\triangle)$. It follows from $\hat{v}\in H^1_0(\Omega)$ and \eqref{prelimeq3},
\begin{align*}
	\|\nabla s_{\hat{v}}\| \le C_1 \|\nabla \hat{v}\|.
\end{align*}
This results in
$$\dfrac{\int_\Omega \nabla \cdot s_{\hat{v}} w}{\|\nabla s_{\hat{v}}\|}=\dfrac{\int_\Omega \nabla \cdot \hat{v} w}{\|\nabla s_{\hat{v}}\|}\geq \dfrac{\int_\Omega \nabla \cdot \hat{v} w}{\sqrt{C_1} \|\nabla \hat{v}\|}\geq C_\beta \|w\|$$
for a positive constant $C_\beta.$
\begin{theorem}\label{convgeup}
	Assume $\mu>0$ and $\textbf{u}-\textbf{u}_s=0$ on $\partial \Omega$. Then, the following inequalities hold:
	\begin{align*}
		\|\nabla (u_s-u)\|
		&\le  C|\triangle|(|u|_{2,2,\Omega} +\frac{1}{\mu}|p|_{1,2,\Omega})\\
		\|p_s-p\| &\le C|\triangle|(\mu |u|_{2,2,\Omega}+|p|_{1,2,\Omega}).
	\end{align*}
\end{theorem}
\begin{proof}
	Assume $\textbf{e}=0$ on $\partial \Omega$. Applying Green's theorem, H\"older's inequality, the divergence theorem, and Lemma \ref{lem1} yields
	\begin{align*}
		\int_\Omega \nabla \cdot v (p_s-s_p) &=\mu \int_\Omega \nabla (u_s-s_u) \nabla v-\mu \int_\Omega \nabla (u-s_u) \nabla v+\int_\Omega \nabla \cdot v (p-s_p)\\
		&\le \mu \|\nabla (u_s-s_u)\| \|\nabla v\|+\mu \|\nabla (u-s_u)\| \|\nabla v\|+ \|p-s_p\| \|\nabla v\|\\
		&\le \|\nabla v\| (\mu \|\nabla (u_s-s_u)\| +\mu \|\nabla (u-s_u)\|+ \|p-s_p\|)\\
		&\le \|\nabla v\| (\mu \|\nabla (u_s-s_u)\| +C|\triangle|(\mu |u|_{2,2,\Omega}+|p|_{1,2,\Omega})), \forall v\in H^1_0(\Omega)
	\end{align*}
	which implies by \eqref{infsup_inequpre}
	\begin{align}\label{ineqp}
		\|p_s-s_p\| \le  \mu \|\nabla (u_s-s_u)\| +C|\triangle|(\mu |u|_{2,2,\Omega}+|p|_{1,2,\Omega}).
	\end{align}
	Furthermore, we have
	\begin{align*}
		&\mu \|\nabla (u_s-s_u)\|^2=\mu \int_\Omega \nabla (u_s-s_u) \nabla (u_s-s_u)\\&= \mu \int_\Omega \nabla (u-s_u) \nabla (u_s-s_u)- \int_\Omega \nabla \cdot (u_s-s_u)  (p-s_p)+\int_\Omega \nabla \cdot (u_s-s_u)  (p_s-s_p)\\
		&= \mu \int_\Omega \nabla (u-s_u) \nabla (u_s-s_u)- \int_\Omega \nabla \cdot (u_s-s_u)  (p-s_p)+\int_\Omega \nabla \cdot (u-s_u)  (p_s-s_p)\\
		&\le \mu \|\nabla (u-s_u)\|\|\nabla (u_s-s_u)\|+\|\nabla (u_s-s_u)\|\|p-s_p\|+\|\nabla (u-s_u)\|\|p_s-s_p\| \\
		&\le \|\nabla (u_s-s_u)\|(\mu \|\nabla (u-s_u)\|+\|p-s_p\|)\\&+\|\nabla (u-s_u)\|( \mu \|\nabla (u_s-s_u)\| +C|\triangle|(\mu |u|_{2,2,\Omega}+|p|_{1,2,\Omega})) 
	\end{align*}
	We obtain
	\begin{align*}
		\mu \|\nabla (u_s-s_u)\|^2\le \|\nabla (u_s-s_u)\|(2 \mu \|\nabla (u-s_u)\|+\|p-s_p\|)+\|\nabla (u-s_u)\|C|\triangle|(\mu |u|_{2,2,\Omega}+|p|_{1,2,\Omega})
	\end{align*}
	By simple calculation, we get $\mu A^2-BA-C\le 0$, then $A\le \frac{1}{\mu} B+\sqrt{\frac{1}{\mu} C}$,
	Then, we have 
	\begin{align*}
		\|\nabla (u_s-s_u)\|\le 2\|\nabla (u-s_u)\|+\frac{1}{\mu} \|p-s_p\|+\sqrt{\|\nabla (u-s_u)\|C|\triangle||u|_{2,2,\Omega}+\frac{1}{\mu} \|\nabla (u-s_u)\|C|\triangle| |p|_{1,2,\Omega}}
	\end{align*}
	and by Lemma \ref{lem1} and Young's inequality we get
	\begin{align*}
		\|\nabla (u_s-s_u)\|&\le 2C|\triangle||u|_{2,2,\Omega}+\frac{1}{\mu} C|\triangle||p|_{1,2,\Omega}+\sqrt{C^2|\triangle|^2|u|^2_{2,2,\Omega}+\frac{1}{\mu} C^2|\triangle|^2|u|_{2,2,\Omega} |p|_{1,2,\Omega}}\\
		&\le C_1|\triangle||u|_{2,2,\Omega}+\frac{1}{\mu} C_1|\triangle||p|_{1,2,\Omega}.
	\end{align*}
	Therefore, we can conclude that
	\begin{align*}
		\|\nabla (u_s-s_u)\|
		&\le  C_1|\triangle||u|_{2,2,\Omega} +\frac{C_1}{\mu}|\triangle||p|_{1,2,\Omega}
	\end{align*}
	for some constant $C_1.$
	By the above inequality and \eqref{ineqp}, we get
	\begin{align*}
		\|p_s-s_p\| &\le  \mu (C|\triangle||u|_{2,2,\Omega} +\frac{C}{\mu}|\triangle||p|_{1,2,\Omega}) +C|\triangle|(\mu |u|_{2,2,\Omega}+|p|_{1,2,\Omega})\\
		&=C_1|\triangle|(\mu |u|_{2,2,\Omega}+|p|_{1,2,\Omega})
	\end{align*}
	for some constant $C_1.$ 
	By Lemma \ref{lem1}, we can finish the proof.
\end{proof}

Assume our domain $\Omega$ is bounded and connected with a Lipschitz continuous boundary. Let us define a new norm $\|\textbf{u}\|_{L}$ on $\tilde{H}_0$ as follows:
\begin{equation}
	\label{H2norm3}
	\| (\textbf{u}, p ) \|_{L}=\| \mu \Delta \textbf{u}+\nabla p\|_{L^2(\Omega)}+\| \nabla \cdot \textbf{u} \|_{L^2(\Omega)}
\end{equation}
We can show that $\|\cdot \|_{L}$ is a norm on $\tilde{H}_0$ as follows:  \\
Indeed, it is clear that $\|a(\textbf{u}, p ) \|_{L}=|a| \|(\textbf{u}, p ) \|_{L}$ and $\|(u_1, u_2, p ) +(v_1, v_2, q ) \|_{L}\le \|(u_1, u_2, p ) \|_{L}+\|(v_1, v_2, q ) \|_{L}$ by the linearity of $\Delta, \nabla, \nabla \cdot$. 
If $\|(u_1, u_2, p ) \|_{L}=0$, then $ \mu \Delta \textbf{u}+\nabla p= 0$ and $\nabla \cdot \textbf{u}=0$ in $\Omega$ and $ \textbf{u}=0$ on the boundary $\partial \Omega$. 
By theorem 5.1 in \cite{GR86}, $\textbf{u}, p$ are zero function. Therefore, $\|\cdot \|_{L}$ is a norm on $\tilde{H}_0$.
We mainly show that the above two norms are equivalent on $\tilde{H}_0$. 
\begin{lemma}
	\label{mjlai2}
	Let $\Omega$ be the connected and bounded domain. Then, there exist two positive constants $A$ and $B$ such that 
	\begin{equation}
		\label{equivalentstokes1}
		A \|(\textbf{u}, p) \|_{H} \le \|(\textbf{u}, p) \|_{L} \le B \| (\textbf{u}, p)\|_{H}, \quad \forall (\textbf{u}, p)\in \tilde{H}_0.
	\end{equation}
\end{lemma}
\begin{proof}
	First, we let $d=2$ and we show that $\|(u_1, u_2, p) \|_{L} \le B \| (u_1, u_2, p)\|_{H}$. 
	We have by simple calculation
	\begin{align*}
		\|(u_1, u_2, p) \|_{L} &\le \|\mu\|_{L^\infty(\Omega)}\| \Delta u_1 \|_{L^2(\Omega)}+\|p_x\|_{L^2(\Omega)}+\|\mu\|_{L^\infty(\Omega)}\| \Delta u_2 \|_{L^2(\Omega)}+\|p_y\|_{L^2(\Omega)}+\| \nabla \cdot \textbf{u} \|_{L^2(\Omega)}\\
		&\le C(1+ \|\mu\|_{L^\infty(\Omega)})(\sum_{0\le |\beta| \le 2} \|D^\beta u_1\|_{L^2(\Omega)}+\sum_{0\le |\beta| \le 2} \|D^\beta u_2\|_{L^2(\Omega)})+C\|p\|_{H^1(\Omega)}\\&\le C_1\|(u_1, u_2, p)\|_{H}.
	\end{align*}
	for some constant $C$.
	By Lemma 5 in \cite{LL22}  and the above inequality, there exist $\alpha >0$ satisfying 
	$$   
	\|(u_1, u_2, p)\|_{H} \le \alpha \|(u_1, u_2, p)\|_L. 
	$$
	Therefore, we choose $A=\frac{1}{\alpha}$ to finish the proof. Similarly, we can proof the theorem when $d=3.$
\end{proof}
\begin{theorem}
	\label{lem3}
	Suppose that $u_1, u_2 \in H^3(\Omega), \textbf{u}-\textbf{u}_{s}=0$ and $p \in H^2(\Omega)$. Then there exists a positive constant $\hat{C}$ depending on $D\ge k$ such that 
	$$
	\| (u_1, u_2, p)-(u_{s_1}, u_{s_2}, u_{s_p})\|_{L}< \hat{C}_{|\triangle|}
	$$
	where $\hat{C}_{|\triangle|}$ is a constant depending on $|\triangle|$.
\end{theorem}
\begin{proof}
	Indeed, we have $s_{u}, s_{p}$ satisfying 
	\begin{align}\label{3ineq1}
		\| -\mu\Delta u_1+p_x+\mu \Delta s_{u_1}- (s_{p})_x\|\leq 
		\frac{C}{2}|\triangle|(\|\mu\|_{L^\infty(\Omega)}|u_1|_{3,2,\Omega}+|p|_{2,2,\Omega})
	\end{align}
	by Lemma~\ref{lem1}. Similarly, we have
	\begin{align}\label{3ineq2}
		\| -\mu\Delta u_2+p_y+\mu \Delta s_{u_2}- (s_{p})_y\|\leq 
		\frac{C}{2}|\triangle|(\|\mu\|_{L^\infty(\Omega)}|u_1|_{3,2,\Omega}+|p|_{2,2,\Omega})
	\end{align}
	Since $(u_{s_1}, u_{s_2}, u_{s_p})$ is the solution for (\ref{PDE2stokes}), we have  
	$$
	|-\mu\Delta u_1(\xi_i)+p_x(\xi_i)+\mu \Delta u_{s_1}(\xi_i)- (u_{s_p})_x(\xi_i)|=0
	$$ 
	for any domain points $\{\xi_i\}$ which construct the matrix $K$.
	Note that $-\mu \Delta (\textbf{u}_{s}-\textbf{s}_{u})+\nabla (p_{s}-\textbf{s}_{p})$ is a polynomial over each triangle $\triangle$ which has small values at the domain points. That is, it is small over $\triangle.$
	This  implies that 
	\begin{align}
		\label{6}
		\|-\mu \Delta (\textbf{u}_{s}-\textbf{s}_{u})+\nabla (p_{s}-\textbf{s}_{p})\| \le \epsilon_\triangle |\Omega|
	\end{align}
	by using Theorem 2.27 in \cite{LS07}. 
	Finally, we can use \eqref{3ineq1} and \eqref{6} to prove 
	\begin{align*}
		&\| -\mu\Delta u_1+p_x+\mu \Delta u_{s1}- (p_s)_x\|\\
		&\le \| -\mu\Delta u_1+p_x+\mu \Delta s_{u1}- (s_p)_x\|+\| -\mu \Delta s_{u1}+(s_p)_x+\mu \Delta u_{s1}- (p_s)_x\|\\
		&\le C|\triangle|(\|\mu\|_{L^\infty(\Omega)}|u_1|_{3,2,\Omega}+|p|_{2,2,\Omega})+ \epsilon_\triangle |\Omega|
	\end{align*}
	for small constant $\epsilon_\triangle.$
\end{proof}

Using Theorem \ref{mjlai2}, we immediately obtain the following theorem
\begin{theorem}
	\label{mainthm7}
	Suppose $f$ and $g$ are continuous over bounded and connected with a Lipschitz continuous boundary in $\mathbb{R}^d$ for $d=2,3$. Suppose that $(\textbf{u}-\textbf{u}_{s}, p-p_s)\in \tilde{H}^3_0$. Then, we have the following inequality
	\begin{align*}
		\|(\textbf{u}-\textbf{u}_s, p-p_s)\|_{H^2} \le \frac{\hat{C}_{|\triangle|}}{A} 
	\end{align*}
	for a positive constant $C$ depending on $\Omega, D $ and $A$ where $A$ is one of the constants in Theorem~\ref{mjlai2}. 
\end{theorem}

\subsection{Existence Results of NS equation}
We introduce a discrete divergent free subspace $S_h$ of $S_D^r$ as follows:
$$ S_h=\{u \in S_D^r: u=0 ~\text{on}~ \partial \Omega ~\text{and} ~\nabla \cdot u=0~ \text{in}~\Omega \}.$$
Let $F: S_{h} \rightarrow S_{h}$ be a nonlinear map so that for each $w \in S_{h}$, $u_{s} = F(w) \in S_{h}$ is given as solution of the following liner problem:
\begin{align}\label{existNAeq}
	-\mu \Delta u_{s}+w \nabla u_{s}+\nabla p_{s}=f
\end{align}
The map $F$ is continuous and therefore compact in the finite-dimensional space $S_h.$
We will use the following Leray-Schander fixed point theorem 
\begin{lemma}
	Let $F$ be a continuous and compact mapping of a Banach space $X$ into itself, such that the set 
	$$\{x\in X: x=\lambda F(x) ~\text{for some}~ 0\le \lambda \le 1\}$$
	is bounded. Then $F$ has a fixed point.
\end{lemma}
to show the existence of the solution.
Then our formulation can be reformulated as seeking $u_s \in S_h$ such that
$$
\begin{aligned}
	\int_{\Omega} f v=& \int_{\Omega}-\mu \Delta u_{s} v+ (w\cdot \nabla) u_{s} v+\nabla p_{s} v \\
	& \mu \int_{\Omega} \nabla u_{s} \nabla v+\int_{\Omega}(w \cdot \nabla) u_{s} v-\int_{\Omega} p (\nabla \cdot v)\\
	& =\mu \int_{\Omega} \nabla u_{s} \nabla v+\int_{\Omega}(w\cdot \nabla)  u_{s} v,~ \forall v\in S_h.
\end{aligned}
$$
If we choose $v=u_{\text {s }}$, for any $w\in S_h$, we have
$$
\begin{aligned}
	& \int_{\Omega} f \cdot u_{s}=\mu \int_{\Omega} \nabla u_{s} \nabla u_{s}+\int_{\Omega} w_{1}\left(u_{1 x}^{s} u_{1}^{s}+u_{2 x}^{s} u_{2}^{s}\right)+w_{2}\left(u_{1 y}^{s} u_{1}^{s}+u_{2 y}^{s} u_{2}^{s}\right) \\
	& =\mu \int_{\Omega} \nabla u_{s} \nabla u_{s}+\int_{\Omega} w_1 \left(\left(u_{1}^{s}\right)^{2}+\left(u_{2}^{s}\right)^{2}\right)_{x} +w_{2}\left(\left(u_{1}^{s}\right)^{2}+\left(u_{2}^{s}\right)^{2}\right)_{y} \\
	& =\mu \int_{\Omega} \nabla u_{s} \nabla u_{s}-\int_{\Omega} \operatorname{div} w \cdot\left(\left(u_{1}^{s}\right)^{2}+\left(u_{2}^{s}\right)^{2}\right)=\lambda \mu\|\nabla u_s\|^2_{L^{2}(\Omega)}
\end{aligned}
$$

If $\lambda>0$ and $w$ satisfies $\lambda F(w)=w,$ then
\begin{align*}
	\lambda \mu\|\nabla u_s\|^2_{L^{2}(\Omega)}&=\int_\Omega f  u_s 
\end{align*}
By introducing a mesh-dependent norm
$$\|f\|_{*,h}=\sum_{v\in S_h} \dfrac{\int_\Omega f v}{\|\nabla v\|_{L^{2}(\Omega)}},$$ we obtain 
$\lambda \leq \frac{\|f\|_{*,h}}{\mu\|\nabla u_s\|_{L^{2}(\Omega)}}$.

Thus, $\lambda<1$ holds true for any $u_s$ being on the boundary of the ball $B_{S_h}^\rho:=\{u_s\in S_h: \|\nabla u_s\|<\rho \}$ in $S_{h}$ centered at the origin with radius $\rho>\frac{\|f\|_{*,h}}{\mu}$. Consequently, the Leray-Schauder fixed point theorem implies that the nonlinear map $F$ defined by \eqref{existNAeq} has a fixed point $\mathbf{u}_{s}$ such that,
$$
F\left(\mathbf{u}_{s}\right)=\mathbf{u}_{s},
$$
in $B_{S_h}^\rho$ with radius $\rho>\frac{\|f\|_{*,h}}{\mu}$. The fixed point $\mathbf{u}_{h}$ also is a solution to our problem, which in turn provides a solution to the LL method. 
To show the uniqueness of the solution of \eqref{PDE2stokes}. Let $\textbf{u}_s$ and $\hat{\textbf{u}}_s$ be two solutions of the equation \eqref{PDE2stokes} in $\mathcal{S}^r_D$. Then $\textbf{u}_s(\xi_i)-\hat{\textbf{u}}_s(\xi_i)=0$ for all $i=1,\cdots, N.$ By using Theorem 2.27 in \cite{LS07}, $\textbf{u}_s=\hat{\textbf{u}}_s.$
This can be summarized in the following theorem.

\begin{theorem}
	The system has at most one solution $\mathbf{u}_{s} \in S_{h}$. Moreover, all the solutions satisfy the following estimates:
	$$
	\left\|\nabla \mathbf{u}_{s}\right\| \leq\frac{\|f\|_{*,h}}{\mu}.
	$$
\end{theorem}

\section{Numerical Results of Stokes and Navier-Stokes Equation}
In this section, we employ a variety of triangulations spanning several bounded 2D and 3D domains to evaluate multiple solutions of the Stokes and Navier-Stokes equation, aiming to determine the precision capabilities of the collocation method. We define the norms to measure the errors as:
\begin{align*}
	\begin{cases}
		|u|_{l_2}&=\sqrt{\frac{\sum_{i=1}^{NI} (u(i))^2}{NI} }\cr
		|u|_{h_1}&=\sqrt{\frac{\sum_{i=1}^{NI} (u_x(i))^2+(u_y(i))^2}{NI} }
	\end{cases}
\end{align*}
In these expressions, $u(i) = u(\eta_i)$, $u_x(i) = u_x(\eta_i)$, and $u_y(i) = u_y(\eta_i)$, where $u, u_x, u_y$ represent specific functions. We calculate the errors over $NI = 1501 \times 1501$ uniformly distributed points in the two-dimensional cases and $NI = 101 \times 101 \times 101$ uniformly distributed points in the three-dimensional context. These points, represented as ${(\eta_i)}_{i=1}^{NI}$, are located within the designated domains. The ensuing tables display the $l_2$ and $h_1$ errors, offering a comparison between the approximated spline solutions and the exact solutions.

In our computational experiments, we use multiple CPUs in the computer so that multiple operations can be done simultaneously. For all simulations, we use 10 processors on a parallel computer equipped with a 12th Gen Intel(R) Core(TM) i7-12650H processor running at 2.30 GHz and 16.0 GB of installed RAM.

Furthermore, we consider both continuous and discontinuous viscosity functions in $\mathbb{R}^2$ as delineated in \cite{JKN15}, \cite{WM21}:
\begin{itemize}
	\item $\mu_1  = \nu_{\min}+(\nu_{\max}-\nu_{\min})xy$ where $ \nu_{\min}=1e-06, \nu_{\max}=1$
	\item $\mu_2= \nu_{\min}+(\nu_{\max}-\nu_{\min})\frac{721}{16}x^2y^2(1-x)(1-y)$ where $ \nu_{\min}=1e-06, \nu_{\max}=1$
	\item $\mu_3=\nu_{min}+(\nu_{max}-\nu_{min})e^{-10^{13} ((x-0.5)^{10}+(y-0.5)^{10})}$ where $ \nu_{\min}=1e-06, \nu_{\max}=1$
	\item $\mu_4=\nu_{min}+(\nu_{max}-\nu_{min})(1-e^{-10^{13} ((x-0.5)^{10}+(y-0.5)^{10})})$ where $ \nu_{\min}=1e-06, \nu_{\max}=1$
	\item $\mu_5= \nu_{\min}+(\nu_{\max}-\nu_{\min})16x(1-x)y(1-y) (0.5+\arctan(\dfrac{\kappa(r-(x-0.5)^2-(y-0.5)^2)}{\pi}))$ where $k=2000, r=0.1, \nu_{\min}=0.1, \nu_{\max}=1$ 
	\item $\mu_6= \nu_{\min}+(\nu_{\max}-\nu_{\min})16x(1-x)y(1-y) (0.5+\arctan(\dfrac{\kappa(r-(x-0.5)^2-(y-0.5)^2)}{\pi}))$ where $k=2000, r=0.01, \nu_{\min}=0.1, \nu_{\max}=1.$
	\item $\mu_{7} = 
	\begin{cases} 
		1, & \text{if } x \leq 0.5, \\
		1000, & \text{if } x > 0.5 
	\end{cases}$
	\item $\mu_{8} = 
	\begin{cases} 
		1, & \text{if } x \leq 0.5, \\
		10^{-5}, & \text{if } x > 0.5 
	\end{cases}$
	\item $\mu_{9} = 
	\begin{cases} 
		10^{-6}, & \text{if not within any specified circles}, \\
		1, & \text{if within any of the circles}.\\
	\end{cases}$
\end{itemize}
\begin{figure}[htbp]
	\centering
	\begin{tabular}{ccc}  
		\begin{subfigure}{0.3\linewidth}  
			\includegraphics[width=\linewidth]{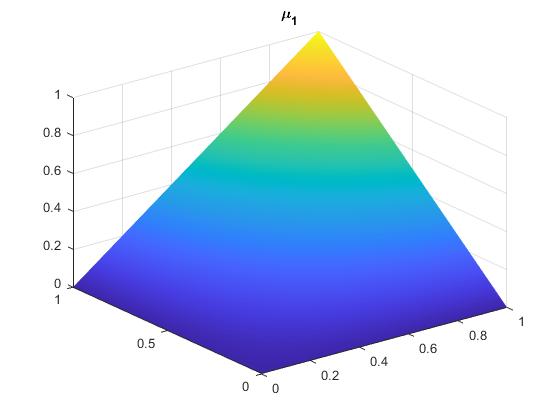}  
			\caption{$\mu_1$}
		\end{subfigure} &
		\begin{subfigure}{0.3\linewidth}
			\includegraphics[width=\linewidth]{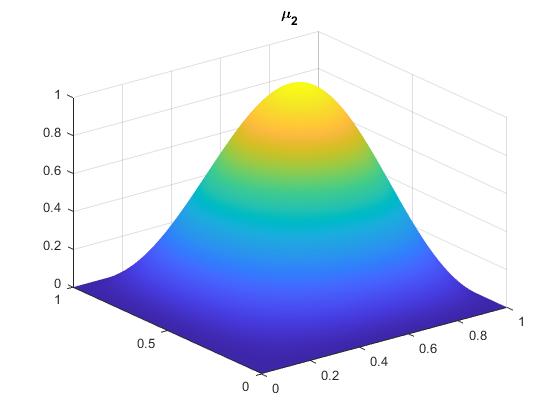}
			\caption{$\mu_2$}
		\end{subfigure} &
		\begin{subfigure}{0.3\linewidth}
			\includegraphics[width=\linewidth]{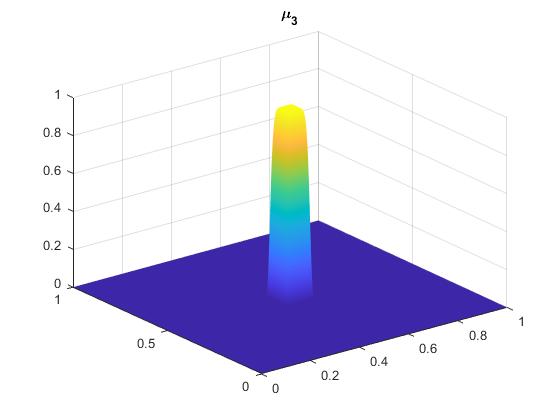}
			\caption{$\mu_3$}
		\end{subfigure} \\
		\begin{subfigure}{0.3\linewidth}
			\includegraphics[width=\linewidth]{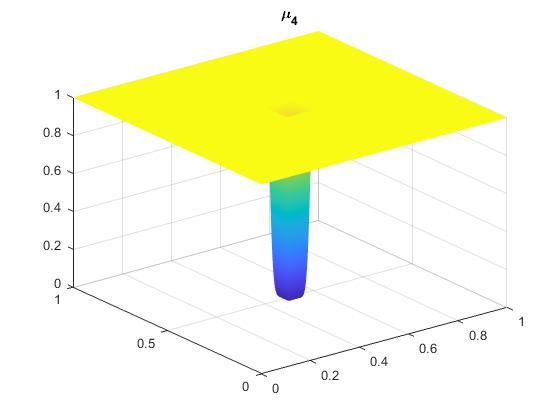}
			\caption{$\mu_4$}
		\end{subfigure} &
		\begin{subfigure}{0.3\linewidth}
			\includegraphics[width=\linewidth]{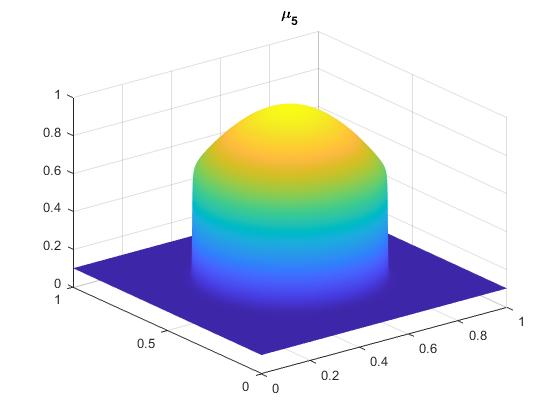}
			\caption{$\mu_5$}
		\end{subfigure} &
		\begin{subfigure}{0.3\linewidth}
			\includegraphics[width=\linewidth]{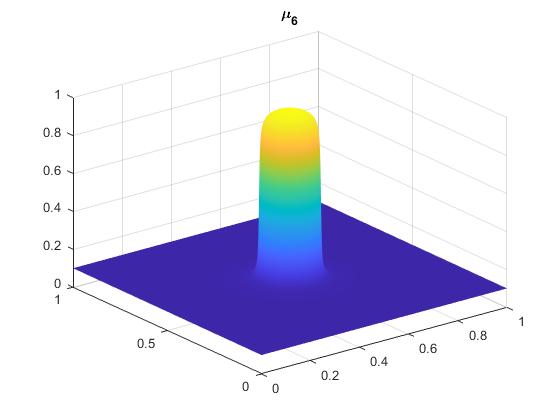}
			\caption{$\mu_6$}
		\end{subfigure} \\
		\begin{subfigure}{0.3\linewidth}
			\includegraphics[width=\linewidth]{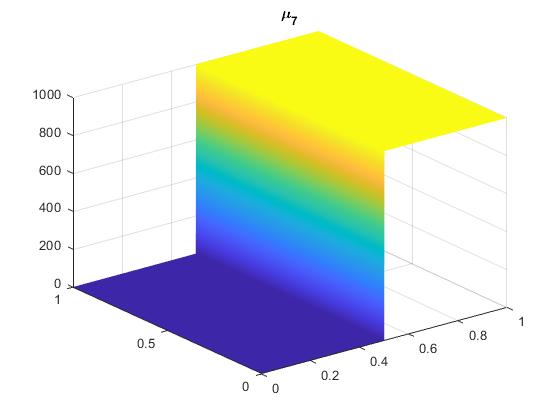}
			\caption{$\mu_7$}
		\end{subfigure} &
		\begin{subfigure}{0.3\linewidth}
			\includegraphics[width=\linewidth]{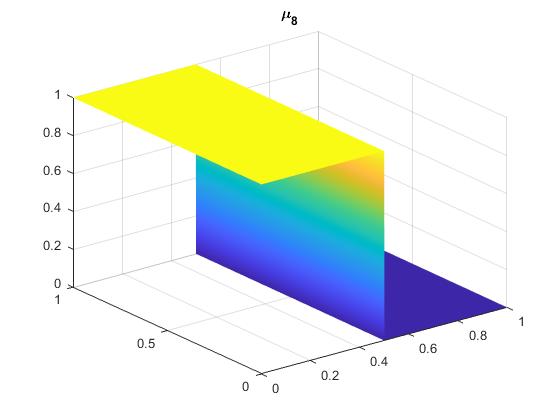}
			\caption{$\mu_8$}
		\end{subfigure} &
		\begin{subfigure}{0.3\linewidth}
			\includegraphics[width=\linewidth]{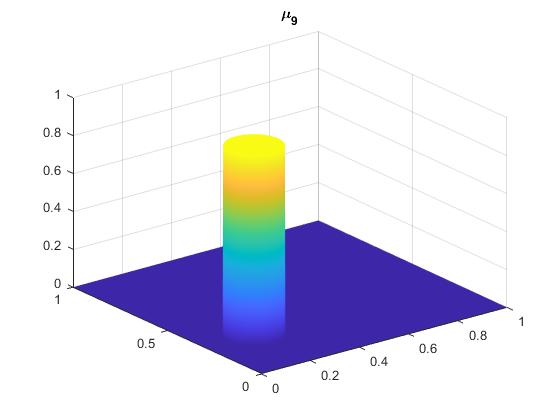}
			\caption{$\mu_9$}
		\end{subfigure}
	\end{tabular}
	\caption{Graph of different $\mu=\mu_1, \cdots, \mu_9$}\label{figure:mu}
\end{figure}
\subsection{Computational experiments on 2D and 3D Stokes Equations}
We conducted a series of computational experiments to evaluate the accuracy of the collocation method in solving the Stokes equations. These experiments involved diverse triangulations across several bounded domains. For brevity, only select computational results are showcased. The tables in this section detail the maximal discrepancies and RMSE when contrasting the spline approximations with the exact solutions.

\begin{example}\label{ex1stokes}
	Consider problem \eqref{stokes} with $\Omega=(0,1)^2.$ The source term and the boundary value g are chosen so that the exact solution is 
	$$\begin{pmatrix} \textbf{u}^{1}(x,y)\\ p^1(x,y)\end{pmatrix}= \begin{pmatrix} 10x^2y(2y - 1)(x - 1)^2(y - 1)\\ -10xy^2(2x - 1)(x - 1)(y - 1)^2 \\ 10(2x-1)(2y-1)\end{pmatrix}.$$ The uniform triangular mesh is used for testing.  The results, obtained using a uniform triangular mesh and denoting mesh size by $h$, are compared with those from the MYZ method \cite{MYZ21}. Table\ref{table:Stokes1} shows the numerical results for each continuous and discontinuous $\mu.$ Table \ref{MYZcomp} shows that we get better results using $D=7, r=2.$ We can observe a trend regarding the errors in our calculations for $\textbf{u}$ and $p$. As $\mu$ gets smaller, the errors in $\textbf{u}$ get larger, but the errors in $p$ get smaller. 
	\begin{figure}[htp] 
		\includegraphics[width=\linewidth]{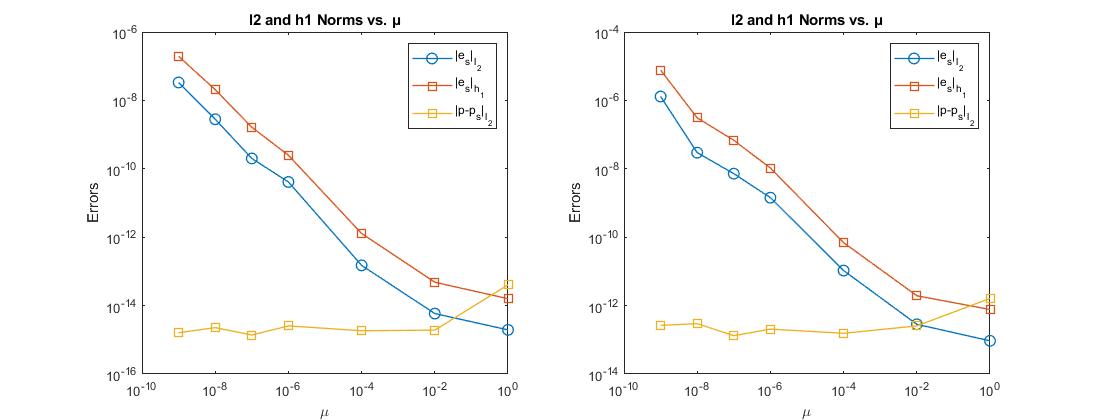}
		\caption{The $l_2, h_1$ errors of $e_s, p-p_s$ for the spline solution of $(\textbf{u}^{1},p^{1})$(Left) and $(\textbf{u}^{2},p^{2})$(Right) for 2D Stokes equations on the domain $[0,1]^2$ for each $\mu=1,1e-02, \cdots, 1e-09$ using the LL method with parameters $D=7, h=1/2$, and $r=2$.}\label{figure:stokeserr}
	\end{figure}
	\begin{table}[htp!]
		\centering
		\small
		\begin{tabular}{ c | c c c c c  } 
			\hline
			$\mu$&CPU time&$\|\nabla \cdot\textbf{u}\|_\infty$&$|e_s|_{l_2} $&$|e_s|_{h_1} $& $|p-p_s|_{l_2} $\\
			\hline  
			$ \mu_1 $& 0.29 &  8.37e-13 & 4.38e-14 & 3.71e-13 & 1.91e-13   \\   
			$ \mu_2 $& 0.28 &  7.49e-13 & 4.73e-14 & 4.42e-13 & 5.26e-14   \\   
			$ \mu_3 $& 0.28 &  6.66e-13 & 5.52e-10 & 4.36e-09 & 2.34e-14   \\   
			$ \mu_4 $& 0.28 &  1.21e-12 & 6.31e-14 & 3.90e-13 & 1.01e-12   \\   
			$ \mu_5 $& 0.29 &  8.85e-13 & 5.31e-14 & 4.81e-13 & 1.28e-13   \\   
			$ \mu_6 $& 0.28 &  1.36e-12 & 5.51e-14 & 3.63e-13 & 8.07e-14   \\   
			$ \mu_7 $& 0.29 &  1.02e-12 & 9.83e-14 & 7.65e-13 & 8.49e-10   \\   
			$ \mu_8 $& 0.28 &  5.03e-13 & 3.54e-12 & 3.83e-11 & 9.59e-13   \\   
			$ \mu_9 $& 0.29 &  5.36e-13 & 4.80e-10 & 3.83e-09 & 2.29e-14   \\    
			\hline
		\end{tabular}
		\caption{The $l_2, h_1$ errors of $e_s, p-p_s$ for the spline solution of $(\textbf{u}^{1},p^{1})$ for 2D Stokes equations on the domain $[0,1]^2$. Computations are carried out using the LL method with parameters $D=7, h=1/2$, and $r=2$ for each different $\mu=\mu_1,\cdots, \mu_9$. The table also presents the CPU time and the divergence norm of $\textbf{u}_s$.}\label{table:Stokes1}
	\end{table}
	\begin{table}
		\centering
		\footnotesize
		\begin{tabular}{ c | c c c c c |c c c } 
			\hline
			&\multicolumn{5}{c|}{LL method}&\multicolumn{3}{c}{MYZ method}\cr
			\hline
			$\mu$&CPU time&$\|\nabla \cdot\textbf{u}\|_\infty$&$|e_s|_{l_2} $&$|e_s|_{h_1} $& $|p-p_s|_{l_2} $&$|e_s|_{l_2} $&$|e_s|_{h_1} $& $|p-p_s|_{l_2} $\\
			\hline
			$\mu=1e+00$ &  0.29 & 1.62e-14 & 1.95e-15 & 1.59e-14 & 4.05e-14 & 2.83e-09&2.10e-06&9.68e-06  \cr  
			$\mu=1e-02$ &0.28 & 3.08e-14 & 5.85e-15 & 4.82e-14 & 1.92e-15 &2.83e-09&2.10e-06&9.68e-08  \cr  
			$\mu=1e-04$ &  0.29 & 2.13e-14 & 1.50e-13 & 1.28e-12 & 1.84e-15&2.83e-09&2.06e-06&1.34e-10 \cr  
			$\mu=1e-06$& 0.28 & 3.49e-14 & 4.12e-11 & 2.49e-10 & 2.55e-15 &2.35e-09&2.06e-06&1.34e-12 \cr  
			$\mu = 1e-07 $ & 0.29 & 2.60e-14 & 2.01e-10 & 1.64e-09 & 1.38e-15 &&&\\   
			$\mu =1e-08  $ & 0.28 & 5.99e-14 & 2.81e-09 & 2.06e-08 & 2.27e-15&&&\\   
			$\mu = 1e-09 $ &  0.29 & 2.97e-14 & 3.34e-08 & 1.94e-07 & 1.61e-15 &&&\\ 
			\hline
		\end{tabular}
		\caption{Comparison of the \(l_2\) and \(h_1\) error measures for \(e_s\) and \(p-p_s\) of the spline solution for \((\mathbf{u}^{1},p^{1})\) related to Stokes equations over the domain \([0,1]^2\). The LL method is implemented with parameters \(D=7\), \(h=\frac{1}{2}\), and \(r=2\), while the MYZ method from \cite{MYZ21} uses \(h=\frac{1}{64}\).
		}\label{MYZcomp}
	\end{table}
\end{example}
\begin{example}
	Let $\Omega=[0,1]^2$, and the exact solutions are chosen as follows:
	$$\begin{pmatrix} \textbf{u}^{2}(x,y)\\ p^2(x,y)\end{pmatrix}= \begin{pmatrix} 0\\0 \\ -\frac{Ra}{2}y^2+Ra y-\frac{Ra}{3},\end{pmatrix}$$
	where $Ra=1000$. Dirichlet boundary condition has been chosen for the test. Table \ref{table:Stokes2} shows the numerical results for several $\mu$. The results in Table \ref{table:Stokes2} underscore the speed and precision of our method.
	\begin{table}[htp!]
		\centering
		\small
		\begin{tabular}{ c | c c c c c  } 
			\hline
			$\mu$&CPU time&$\|\nabla \cdot\textbf{u}\|_\infty$&$|e_s|_{l_2} $&$|e_s|_{h_1} $& $|p-p_s|_{l_2} $\\
			\hline
			$ \mu_1 $& 0.44 & 2.37e-12 & 1.13e-13 & 9.36e-13 & 4.84e-13   \\   
			$ \mu_2 $& 0.44 & 8.82e-13 & 8.95e-14 & 8.43e-13 & 2.96e-13   \\   
			$ \mu_3 $& 0.46 & 1.67e-12 & 1.38e-09 & 8.20e-09 & 1.69e-13   \\   
			$ \mu_4 $& 0.45 & 5.76e-13 & 9.07e-14 & 7.25e-13 & 1.86e-12   \\   
			$ \mu_5 $& 0.44 & 1.20e-12 & 8.31e-14 & 7.55e-13 & 3.20e-13   \\   
			$ \mu_6 $& 0.44 & 1.40e-12 & 1.08e-13 & 8.39e-13 & 3.15e-13   \\   
			$ \mu_7 $& 0.45 & 4.79e-13 & 7.79e-14 & 6.83e-13 & 9.95e-10   \\   
			$ \mu_8 $& 0.47 & 5.72e-13 & 7.32e-12 & 8.77e-11 & 6.77e-13   \\   
			$ \mu_9 $& 0.45 & 6.69e-13 & 3.93e-10 & 4.03e-09 & 3.39e-13   \\   
			\hline
		\end{tabular}
		\caption{The $l_2, h_1$ errors of $e_s, p-p_s$ of spline solution for $(\textbf{u}^{2},p^{2})$ for Stokes equations  over $[0,1]^2$  using LL method with $D=7, h=1/2$ and $r=2$}\label{table:Stokes2}
	\end{table}
\end{example}
\begin{example}
	Retaining the source term and boundary value from Example 1, we extend our examination to diverse domains, as depicted in Figure \ref{fig:2dstokesdomain}, utilizing parameters $D=7, r=2$. The errors, cataloged in Table \ref{table:Stokes3}, attest to the robustness and precision of our collocation method across varied domains.
	\begin{figure}[htbp]
		\centering
		\includegraphics[width=\linewidth]{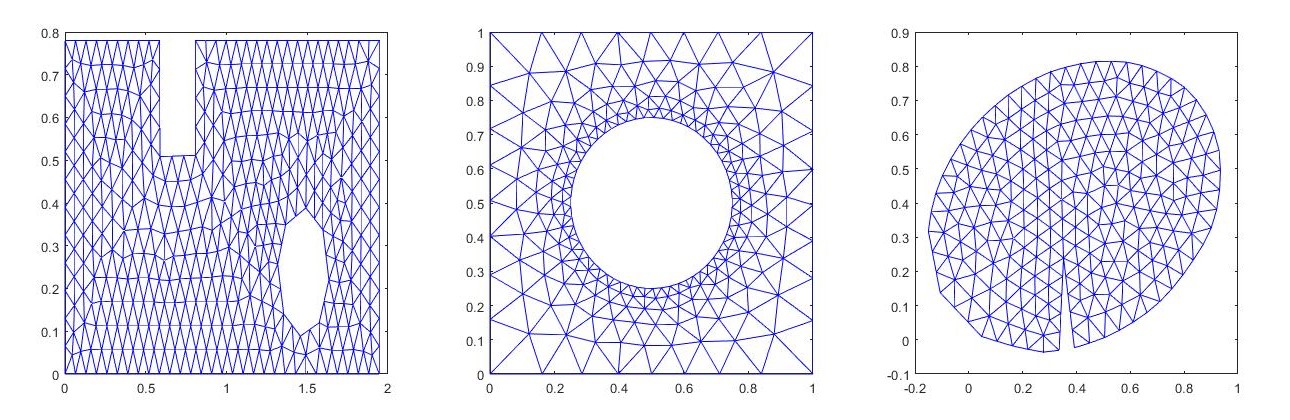}
		\caption{Various 2D domains for the Stokes equation: $D_1$ (left), $D_2$ (middle), and $D_3$ (right).}
		\label{fig:2dstokesdomain}
	\end{figure}
	\begin{table}
		\centering
		\footnotesize
		\begin{tabular}{ c | c c c | c c c | c c c } 
			\hline
			$\mu$ & \multicolumn{3}{c|}{$D_1$} & \multicolumn{3}{c|}{$D_2$} & \multicolumn{3}{c}{$D_3$} \\
			
			\hline
			&$|e_s|_{l_2} $&$|e_s|_{h_1} $& $|p-p_s|_{l_2} $&$|e_s|_{l_2} $&$|e_s|_{h_1} $& $|p-p_s|_{l_2} $&$|e_s|_{l_2} $&$|e_s|_{h_1} $& $|p-p_s|_{l_2} $\\
			\hline
 1e+00& 1.25e-08 & 4.83e-07 & 1.21e-06& 7.90e-07 & 2.17e-05 & 5.43e-05& 4.90e-08 & 1.30e-06 & 3.68e-06   \\   
1e-02& 5.05e-09 & 1.46e-07 & 4.71e-09& 2.16e-07 & 6.97e-06 & 1.01e-07& 1.91e-08 & 3.97e-07 & 2.23e-08   \\   
1e-04& 2.12e-07 & 4.91e-06 & 1.20e-09& 1.75e-06 & 5.03e-05 & 6.80e-09& 8.06e-07 & 1.81e-05 & 5.99e-09   \\   
1e-06& 1.73e-05 & 2.75e-04 & 1.05e-10& 1.53e-05 & 2.93e-04 & 3.64e-10& 5.32e-05 & 8.44e-04 & 1.31e-09   \\   
1e-07& 5.85e-05 & 7.93e-04 & 2.02e-11& 5.45e-05 & 8.54e-04 & 1.58e-10& 2.83e-04 & 4.17e-03 & 8.15e-10   \\   
1e-08& 6.19e-05 & 8.30e-04 & 2.08e-12& 2.49e-04 & 3.05e-03 & 4.29e-11& 4.67e-04 & 6.70e-03 & 1.46e-10   \\   
1e-09& 6.19e-05 & 8.31e-04 & 2.08e-13& 2.71e-04 & 3.31e-03 & 4.65e-12& 4.81e-04 & 6.88e-03 & 1.48e-11   \\   
$ \mu_1 $& 6.19e-09 & 2.43e-07 & 2.81e-07& 5.09e-07 & 1.24e-05 & 1.85e-05& 7.43e-08 & 1.11e-06 & 2.74e-06   \\   
$ \mu_2 $& 8.43e-09 & 2.71e-07 & 3.29e-07& 6.14e-07 & 1.54e-05 & 1.70e-06& 4.84e-07 & 1.73e-05 & 1.15e-03   \\   
$ \mu_3 $& 1.61e-05 & 2.46e-04 & 4.15e-10& 1.53e-05 & 2.93e-04 & 3.64e-10& 4.88e-05 & 8.11e-04 & 1.71e-09   \\   
$ \mu_4 $& 1.25e-08 & 4.82e-07 & 1.21e-06& 7.90e-07 & 2.17e-05 & 5.43e-05& 4.86e-08 & 1.30e-06 & 3.71e-06   \\   
$ \mu_5 $& 9.12e-09 & 2.20e-07 & 2.53e-07& 2.22e-07 & 7.95e-06 & 2.09e-06& 1.78e-08 & 3.94e-07 & 4.61e-07   \\   
$ \mu_6 $& 5.00e-09 & 1.69e-07 & 5.12e-08& 2.29e-07 & 7.97e-06 & 2.09e-06& 1.21e-08 & 3.52e-07 & 1.24e-07   \\   
$ \mu_7 $& 4.27e-06 & 7.81e-05 & 6.23e-01& 4.42e-05 & 4.03e-04 & 1.78e+00& 3.68e-05 & 3.49e-04 & 6.25e+00   \\   
$ \mu_8 $& 1.54e-06 & 3.06e-05 & 9.94e-07& 2.46e-05 & 4.69e-04 & 4.11e-05& 7.37e-06 & 1.54e-04 & 9.11e-07   \\   
$ \mu_9 $& 1.05e-04 & 1.27e-03 & 1.12e-07& 1.52e-04 & 2.90e-03 & 6.27e-08& 8.81e-05 & 1.30e-03 & 1.49e-07   \\   
			\hline
		\end{tabular}
		\caption{The $l_2, h_1$ errors of $e_s, p-p_s$ of spline solution for $(\textbf{u}^{1},p^{1})$ for Stokes equations  over several domains in Figure \ref{fig:2dstokesdomain}  using LL method with $D=5$ and $r=2$}\label{table:Stokes3}
	\end{table}
\end{example}
\begin{example}
	We explore the domain $\Omega=[-1/2,3/2]\times [0,2]$ with the exact solution
	$$\begin{pmatrix} \textbf{u}^{3}(x,y)\\ p^3(x,y)\end{pmatrix}= \begin{pmatrix}  1 - \exp(\lambda x) \cos(2\pi y)\\\frac{\lambda}{2\pi} \exp(\lambda  x) \sin(2\pi y) \\ \frac{1}{2} \exp(2\lambda x)\end{pmatrix}$$
	where $\lambda = \frac{\text{Re}}{2} - \sqrt{\left(\frac{\text{Re}}{2}\right)^2 + 4\pi^2}$ and $\text{Re} = \frac{1}{\mu}$ denotes the Reynolds number. Table \ref{table:exampleRE} delineates the error profiles for various Reynolds numbers.
	\begin{table}
		\centering
		\small
		\begin{tabular}{ c | c c c c c  } 
			\hline
			Re&CPU time&$\|\nabla \cdot\textbf{u}\|_\infty$&$|e_s|_{l_2} $&$|e_s|_{h_1} $& $|p-p_s|_{l_2} $\\
			\hline
 1e+00& 4.44 &  1.28e-02 & 4.76e-05 & 1.01e-03 & 2.36e-03   \\   
1e+02& 4.37 &  4.49e-06 & 6.05e-07 & 3.68e-06 & 6.81e-08   \\   
1e+04& 4.36 &  4.80e-08 & 8.39e-07 & 4.41e-06 & 1.05e-09   \\   
1e+06& 4.53 &  4.80e-08 & 8.55e-07 & 4.55e-06 & 1.11e-11   \\   
1e+07& 4.29 &  4.80e-08 & 8.55e-07 & 4.55e-06 & 1.11e-12   \\   
1e+08& 4.34 &  4.80e-08 & 8.56e-07 & 4.56e-06 & 1.11e-13   \\   
1e+09& 4.22 &  4.80e-08 & 9.00e-07 & 4.64e-06 & 1.17e-14   \\      
			\hline
		\end{tabular}
		\caption{The $l_2, h_1$ errors of $e_s, p-p_s$ of spline solution for $(\textbf{u}^{3},p^{3})$ for Stokes equations  over $[-1/2,3/2]\times [0,2]$  using LL method with $D=7, h=1/4$ and $r=2$}\label{table:exampleRE}
	\end{table}
\end{example}
\begin{example}
	We extend our investigation to a three-dimensional setting, focusing on a flow within a unit cube $[0,1]^3$. The velocity and pressure fields are expressed as
	\[
	\begin{pmatrix} 
		\textbf{u}^{3d1}(x,y,z)\\ 
		p^{3d1}(x,y,z)
	\end{pmatrix} = 
	\begin{pmatrix} 
		\sin(\pi x) \cos(\pi y) - \sin(\pi x) \cos(\pi z)\\ 
		\sin(\pi y) \cos(\pi z) - \sin(\pi y) \cos(\pi x)\\ 
		\sin(\pi z) \cos(\pi x) - \sin(\pi z) \cos(\pi y)\\ 
		\sin(\pi x) \sin(\pi y) \sin(\pi z)
	\end{pmatrix}
	\]
	and
	$$
	\begin{pmatrix} 
		\textbf{u}^{3d2}(x,y,z)\\ 
		p^{3d2}(x,y,z)
	\end{pmatrix} = 
	\begin{pmatrix} 
		-\exp(x + 2y + 3z)\\ 
		2\exp(x + 2y + 3z)\\ 
		-\exp(x + 2y + 3z)\\ 
		\exp(x + y + z)
	\end{pmatrix}.
	$$
	The performance, in terms of computational errors and the corresponding CPU times of $(\textbf{u}^{3d1},p^{3d1})$ and $(\textbf{u}^{3d2},p^{3d2})$ for varied values of $\mu$, is detailed in Tables \ref{table:3dstokes} and \ref{table:3dstokes1}. These results underscore the efficacy of our methodology. The patterns observed in the $l_2$ and $h_1$ errors, in addition to the pressure errors, are consistent with the observations made during our two-dimensional experiments.
	\begin{table}[htp!]
		\centering
		\begin{tabular}{c c c c c}
			\hline
			$\mu$ &CPU time&$|e_s|_{l_2} $&$|e_s|_{h_1} $& $|p-p_s|_{l_2}$ \\
			\hline
			1e+00& 12.86 & 2.78e-05 & 8.34e-04 & 1.37e-03   \\   
			1e-01& 12.97 & 1.19e-05 & 4.05e-04 & 5.61e-05   \\   
			1e-02& 12.48 & 1.90e-05 & 4.40e-04 & 1.36e-05   \\   
			1e-03& 12.54 & 1.14e-04 & 1.11e-03 & 1.39e-05   \\   
			1e-04& 12.96 & 5.89e-04 & 4.91e-03 & 1.38e-05   \\   
			1e-05& 12.75 & 5.70e-03 & 4.57e-02 & 1.37e-05   \\   
			1e-06& 12.36 & 5.87e-02 & 4.65e-01 & 1.37e-05   \\   
			\hline
		\end{tabular}
		\caption{The $l_2, h_1$ errors of $e_s, p-p_s$ of spline solution for $(\textbf{u}^{3d1},p^{3d1})$ for 3D Stokes equations  over $[0,1]^3$  using LL method with $D=7,h=1/2$ and $r=1$}\label{table:3dstokes}
	\end{table}
	\begin{table}[htp!]
		\centering
		\begin{tabular}{c c c c c}
			\hline
			$\mu$ &CPU time&$|e_s|_{l_2} $&$|e_s|_{h_1} $& $|p-p_s|_{l_2}$ \\
			\hline
			1e+00& 12.91 & 1.41e-04 & 3.80e-03 & 8.23e-03   \\   
			1e-01& 12.53 & 4.69e-05 & 1.68e-03 & 2.45e-04   \\   
			1e-02& 12.36 & 2.78e-05 & 8.77e-04 & 1.22e-05   \\   
			1e-03& 12.84 & 1.03e-04 & 1.14e-03 & 2.62e-06   \\   
			1e-04& 12.40 & 6.04e-04 & 5.30e-03 & 1.72e-06   \\   
			1e-05& 13.00 & 2.99e-03 & 2.16e-02 & 7.36e-07   \\   
			1e-06& 12.49 & 4.78e-03 & 3.40e-02 & 1.19e-07   \\   
			\hline
		\end{tabular}
		\caption{The $l_2, h_1$ errors of $e_s, p-p_s$ of spline solution for $(\textbf{u}^{3d2},p^{3d2})$ for 3D Stokes equations  over $[0,1]^3$  using LL method with $D=7,h=1/2$ and $r=1$}\label{table:3dstokes1}
	\end{table}
\end{example}
\subsection{Computational experiments on 2D and 3D Navier-Stokes Equations}
This section underscores the adaptability and consistency of the LL method, demonstrating its applicability in deriving accurate solutions to the Navier-Stokes equations. We test our numerical example for different $\mu$ in the previous section. We choose the same $\nu_{\max}$ and we choose $\nu_{\min}=1e-03$ for $\mu_1\cdots, \mu_6$ in this section. The tables encapsulate detailed error metrics, offering a comprehensive insight into the precision and reliability of the LL method across varied conditions and domains.
\begin{example}
	Consider problem \eqref{nstokes} with $\Omega=(0,1)^2.$ The source term and the boundary value g are chosen so that the exact solution is 
	$$\begin{pmatrix} \textbf{u}^{1}(x,y)\\ p^1(x,y)\end{pmatrix}= \begin{pmatrix} 10x^2y(2y - 1)(x - 1)^2(y - 1)\\ -10xy^2(2x - 1)(x - 1)(y - 1)^2 \\ 10(2x-1)(2y-1)\end{pmatrix}.$$  The results, displayed in Table \ref{table:NStokes1}, are obtained using a uniform triangular mesh, showcasing the method's precision.
	\begin{table}[htp!]
		\centering
		\small
		\begin{tabular}{ c | c c c c c  } 
			\hline
			$\mu$&CPU time&$\|\nabla \cdot\textbf{u}\|_\infty$&$|e_s|_{l_2} $&$|e_s|_{h_1} $& $|p-p_s|_{l_2} $\\
			\hline
			1e+00& 1.31 & 3.93e-15 & 2.45e-15 & 2.10e-14 & 4.99e-14   \\   
			1e-01& 1.06 & 5.97e-15 & 3.46e-15 & 2.87e-14 & 5.79e-15   \\   
			1e-02& 0.93 & 1.03e-14 & 1.23e-14 & 7.45e-14 & 3.03e-15   \\   
			1e-03& 1.68 & 9.27e-15 & 4.29e-14 & 2.81e-13 & 3.89e-15   \\   
			$ \mu_1 $& 1.42 & 7.19e-15 & 4.40e-15 & 3.86e-14 & 5.40e-14   \\   
			$ \mu_2 $& 1.59 & 7.32e-15 & 4.82e-15 & 3.94e-14 & 9.54e-15   \\   
			$ \mu_3 $& 1.53 & 1.02e-14 & 3.12e-14 & 2.20e-13 & 2.59e-15   \\   
			$ \mu_4 $& 1.33 & 4.52e-15 & 3.30e-15 & 2.89e-14 & 6.05e-14   \\   
			$ \mu_5 $& 2.32 & 7.63e-15 & 8.18e-15 & 9.03e-14 & 2.94e-15   \\   
			$ \mu_6 $& 1.99 & 7.96e-15 & 3.83e-14 & 2.30e-13 & 3.68e-15   \\   
			$ \mu_7 $& 2.09 & 6.83e-15 & 1.30e-14 & 2.19e-13 & 5.01e-14   \\   
			$ \mu_8 $& 2.09 & 5.14e-15 & 6.33e-15 & 7.43e-14 & 4.99e-14   \\   
			$ \mu_9 $& 2.82 & 8.77e-15 & 3.54e-14 & 2.44e-13 & 3.12e-15   \\   
			\hline
		\end{tabular}
		\caption{The $l_2, h_1$ errors of $e_s, p-p_s$ for the spline solution of $(\textbf{u}^{1},p^{1})$ for Navier-Stokes equations on the domain $[0,1]^2$. Computations are carried out using the LL method with parameters $D=7, h=1/2$, and $r=2$. The table also presents the CPU time and the divergence norm of $\textbf{u}$.}\label{table:NStokes1}
	\end{table}
\end{example}
\begin{example}
	Let $\Omega=[0,1]^2$, and the exact solutions are chosen as follows:
	$$\begin{pmatrix} \textbf{u}^{2}(x,y)\\ p^2(x,y)\end{pmatrix}= \begin{pmatrix} 0\\0 \\ -\frac{Ra}{2}y^2+Ra y-\frac{Ra}{3},\end{pmatrix}$$
	where $Ra=1000$. We employ the Dirichlet boundary condition for this evaluation. The numerical errors and rates, obtained with parameters $D=5, \mu=1, r=2$, are detailed in Table \ref{table:NStokes2}. Additionally, Table \ref{table:NStokes3} affirms the effectiveness of our approach in solving the Navier-Stokes equation for different continuous and discontinuous $\mu$.
	\begin{table}[htp!]
		\centering
		\small
		\begin{tabular}{ c | c c c c c c c } 
			\hline
			$h$&CPU time&$|e_s|_{l_2} $&Rate&$|e_s|_{h_1} $&Rate& $|p-p_s|_{l_2} $&Rate\\
			\hline
			1/2& 0.86 & 7.38e-02&-& 5.17e-01 & - & 1.13e+00 &- \\   
			1/4& 0.76& 2.50e-03& 4.88  & 3.25e-02 & 3.99 & 7.08e-02 & 4.00   \\   
			1/8& 2.58  & 8.61e-05& 4.86  & 1.54e-03 & 4.40 & 3.26e-03 & 4.44   \\   
			1/16& 44.89& 5.22e-06& 4.04  & 8.36e-05 & 4.20 & 1.74e-04 & 4.23   \\
			\hline
		\end{tabular}
		\caption{The $l_2, h_1$ errors of $e_s, p-p_s$ of spline solution for $(\textbf{u}^{2},p^{2})$ for Navier-Stokes equations  over $[0,1]^2$  using LL method with $\mu=1, D=5$ and $r=2$}\label{table:NStokes2}
	\end{table}
	\begin{table}[htp!]
		\centering
		\small
		\begin{tabular}{ c | c c c c c  } 
			\hline
			$\mu$&CPU time&$\|\nabla \cdot\textbf{u}\|_\infty$&$|e_s|_{l_2} $&$|e_s|_{h_1} $& $|p-p_s|_{l_2} $\\
			\hline
			1e+00& 0.64 & 1.55e-13 & 8.23e-14 & 8.94e-13 & 1.32e-12   \\   
			1e-01& 0.62 & 2.32e-13 & 9.10e-14 & 8.50e-13 & 1.67e-13   \\   
			1e-02& 0.63 & 2.87e-13 & 3.80e-13 & 2.30e-12 & 9.72e-14   \\   
			1e-03& 0.82 & 4.05e-13 & 1.08e-12 & 1.11e-11 & 8.79e-14   \\   
			$ \mu_1 $& 0.85 & 2.98e-13 & 1.36e-13 & 1.37e-12 & 1.01e-12   \\   
			$ \mu_2 $& 0.89 & 2.61e-13 & 1.36e-13 & 1.29e-12 & 4.13e-13   \\   
			$ \mu_3 $& 0.84 & 3.60e-13 & 1.87e-12 & 1.10e-11 & 8.75e-14   \\   
			$ \mu_4 $& 1.03 & 1.19e-13 & 7.82e-14 & 7.20e-13 & 1.09e-12   \\   
			$ \mu_5 $& 1.01 & 2.38e-13 & 1.66e-13 & 2.22e-12 & 1.99e-13   \\   
			$ \mu_6 $& 0.84 & 3.64e-13 & 5.45e-13 & 5.92e-12 & 8.17e-14   \\   
			$ \mu_7 $& 0.82 & 2.55e-13 & 2.79e-11 & 3.91e-10 & 1.45e-12   \\   
			$ \mu_8 $& 0.84 & 2.54e-13 & 6.19e-13 & 5.06e-12 & 1.67e-12   \\   
			$ \mu_9 $& 0.82 & 3.43e-13 & 1.40e-12 & 9.62e-12 & 9.85e-14   \\     
			\hline
		\end{tabular}
		\caption{The $l_2, h_1$ errors of $e_s, p-p_s$ for the spline solution of $(\textbf{u}^{2},p^{2})$ for Navier-Stokes equations on the domain $[0,1]^2$. Computations are carried out using the LL method with parameters $D=7, h=1/2$, and $r=2$. The table also presents the CPU time and the divergence norm of $\textbf{u}$.}\label{table:NStokes3}
	\end{table}
\end{example}
\begin{example}
	In this example, we focus on solving the Navier-Stokes equation on $[0,1]^2$ with $\mu=1$ and the exact functions are chosen as follows: 
	\[
	\begin{pmatrix}
		\textbf{u}^{3}(x,y) \\
		p^3(x,y) 
	\end{pmatrix} = 
	\begin{pmatrix}
		0.5  \sin^2(2\pi x)  \sin(2\pi y)  \cos(2\pi y) \\
		-0.5 \sin^2(2\pi y)  \sin(2\pi x)  \cos(2\pi x) \\
		\pi^2  \sin(2\pi x)\cos(2\pi y)
	\end{pmatrix}
	\]
	The computed $l_2, h_1$ errors and convergence rates between the spline approximation and the exact solutions are tabulated in Table \ref{table:NStokesconv}. These $l_2, h_1$ errors between the spline approximation and exact solutions are asymptotically proportional to $\mathcal{O}(h^{4.84}).$
	Also, we compare our numerical errors with the errors using the BMV method in \cite{BMV19} and get better result using the same mesh size. Therefore, we conclude that the our method work very well.
	\begin{table}[htp!]
		\centering
		\small
		\begin{tabular}{ c | c c c c c c| c c } 
			\hline
			$h$  & \multicolumn{6}{c|}{LL Method} & \multicolumn{2}{c}{BMV Method} \\
			\hline
			& $|e_s|_{l_2}$ & Rate & $|e_s|_{h_1}$ & Rate & $|p-p_s|_{l_2}$& Rate & $|e_s|_{h_1}$ & $|p-p_s|_{l_2}$ \\
			\hline
			$1/2$ &3.53e-02& -    & 4.40e-01 & -    & 1.04e+00 & -  & - &- \\   
			$1/4$ &1.36e-03& 4.69 & 2.85e-02 & 3.95 & 8.33e-02& 3.64 & -  & -  \\   
			$1/8$ &3.22e-05& 5.40 & 1.24e-03 & 4.52 & 2.48e-03& 5.07  & 3.70e-01 & 3.89e-01  \\   
			$1/16$&1.22e-06& 4.73 & 2.14e-05 & 5.86 & 4.62e-05 & 5.75 & 9.15e-02 & 8.88e-02  \\   
			\hline
		\end{tabular}
		\caption{A comparison of the $l_2, h_1$ errors and convergence rates for the spline solution of $(\textbf{u}^{3},p^{3})$ on the domain $[0,1]^2$. Results from the LL method with parameters $D=7$, and $r=2$ are presented alongside corresponding results obtained using the BMV method introduced in \cite{BMV19}}
		\label{table:NStokesconv}
	\end{table}
	
\end{example}
\begin{example}
	We consider a driven cavity flow with the Dirichlet boundary conditions
	\begin{align*}
		u=\begin{cases}
			(1,0) ~\text{on}~(0,1)\times \{1\},\\
			(0,0) ~\text{otherwise}.
		\end{cases}
	\end{align*}
	As the exact solution is not available, we depict the velocity streamlines, which are computed utilizing a uniform triangular mesh with a granularity of \(h=1/8\), to analyze the Navier-Stokes problem. Figures \ref{figure:cavity1} and \ref{figure:cavity2} exhibit the graphs of velocity and pressure, respectively, ascertained through the LL method for \(\mu=1\) and \(0.01\).
	
	\begin{figure}[htbp]
		\centering
		\includegraphics[width=\linewidth]{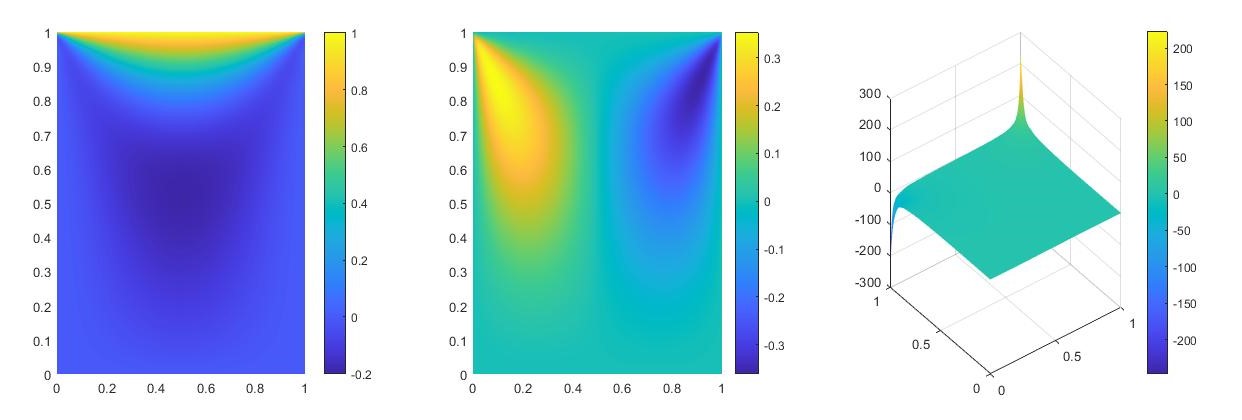}
		\caption{The numerical velocity and pressure when $\mu=1$}\label{figure:cavity1}

	\end{figure}
	\begin{figure}[htbp]
		\centering
		\includegraphics[width=\linewidth]{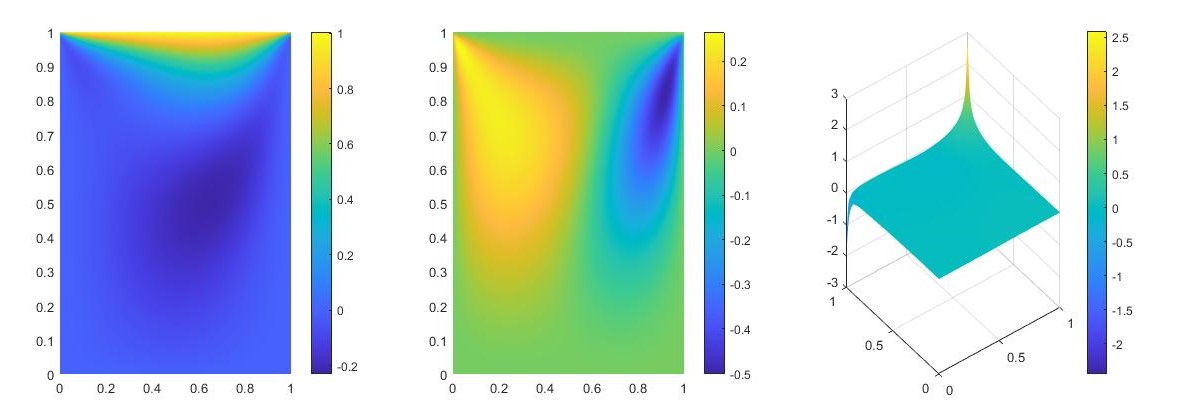}
		\caption{The numerical velocity and pressure when $\mu=0.01$}\label{figure:cavity2}
	\end{figure}
\end{example}
\begin{example}
	Transitioning to a three-dimensional context, we explore the LL method’s performance on the 3D Navier-Stokes equations within the domain $\Omega=[0,1]^3$. 
	We consider the following exact solution
	\[
	\begin{pmatrix} 
		\textbf{u}^{3dns1}(x,y,z)\\ 
		p^{3dns1}(x,y,z)
	\end{pmatrix} = 
	\begin{pmatrix} 
		\sin(z) + \cos(y)\\ 
		\sin(x) + \cos(z)\\ 
		\sin(y) + \cos(x)\\ 
		x
	\end{pmatrix},
	\]
	and 
	\[
	\begin{pmatrix} 
		\textbf{u}^{3dns2}(x,y,z)\\ 
		p^{3dns2}(x,y,z)
	\end{pmatrix} = 
	\begin{pmatrix} 
		-\exp(x + 2y + 3z)\\ 
		2\exp(x + 2y + 3z)\\ 
		-\exp(x + 2y + 3z)\\ 
		\exp(x + y + z)
	\end{pmatrix}.
	\]
	Table \ref{table:NStokes4} presents the numerical errors associated with velocity and pressure for each degree \(D\). Additionally, Table \ref{table:NStokes5} provides a comparative analysis of the \(L^2\) error between the velocity and pressure obtained from the AL method, as cited in \cite{AL02}, and our LL method across various mesh sizes. The results from both tables collectively affirm the efficacy of our LL method.
	
	\begin{table}[htp!]
		\centering
		\small
		\begin{tabular}{ c | c c c } 
			\hline
			Degree&$\|\nabla \cdot\textbf{u}\|_\infty$&$|e_s|_{l_2} $& $|p-p_s|_{l_2} $\\
			\hline
			$D=3$&3.35e-03 & 1.11e-03 & 5.19e-02  \cr  
			$D=4 $&2.45e-04 & 6.20e-05 & 1.96e-03  \cr  
			$D=5 $&1.55e-05 & 2.35e-06 & 1.10e-04  \cr  
			$D=6 $&2.22e-07 & 3.08e-08 & 1.12e-06  \cr  
			\hline
		\end{tabular}
		\caption{The $l_2$ errors and the divergence norm of $e_s, p-p_s$ for the spline solution of $(\textbf{u}^{3dns1},p^{3dns1})$ for Navier-Stokes equations on the domain $[0,1]^3$. Computations are carried out using the LL method with parameters $h=1/2$, and $r=2$ for each degree $D=3,4,5,6$.}\label{table:NStokes4}
	\end{table}
	\begin{table}[htp!]
		\centering
		\small
		\begin{tabular}{c| c  c |c c } 
			\hline
			&\multicolumn{2}{c|}{LL method}&\multicolumn{2}{c}{AL method}\cr
			\hline
			mesh size($h$)&$|e_s|_{l_2} $& $|p-p_s|_{l_2} $&$|e_s|_{l_2} $& $|p-p_s|_{l_2} $\\
			\hline
			1&3.44e-02 & 1.85&3.96e-01&27.40  \cr  
			1/2& 1.92e-03 & 1.04e-01&8.51e-03&7.18e-01  \cr  
			\hline
		\end{tabular}
		\caption{The $l_2$ errors and the divergence norm of $e_s, p-p_s$ for the spline solution of $(\textbf{u}^{3dns2},p^{3dns2})$ for Navier-Stokes equations on the domain $[0,1]^3$. Computations are carried out using the LL method with parameters $h=1, 1/2$, and $r=2$.}\label{table:NStokes5}
	\end{table}
\end{example}

\section{Strategies for Solving Stokes Equations in Domains with Curved Boundaries}

In this section, we explore various techniques for solving the Stokes equation within domains enclosed by curved boundaries. Recently, Isogeometric Analysis (IGA) has gained popularity. It is a method based on utilizing functions that align with those in Computer-Aided Design (CAD) systems to define the domain \( \Omega \). Another notable approach is presented in \cite{S19}, where the Immersed Penalized Boundary Method (IPBM) is employed for solving both two-dimensional and three-dimensional partial differential equations, positioning spline functions as an alternative to IGA.

A significant attribute that distinguishes IPBM from IGA is its non-reliance on mapping between a computational domain and the target domain \( \Omega \). IPBM imposes boundary conditions directly at points located on the boundary \( \partial \Omega \), as opposed to nearby points. This feature differentiates it from many immersion methods discussed in the existing literature, which typically modify or stabilize approximation functions close to \( \partial \Omega \). Additionally, there is no necessity to increase mesh density near the boundary or to compute integrals along it. 

We now elaborate on the application of IPBM in our study. Initially, we consider a mesh that encompasses the curved domain \( \Omega \). We then construct the basis functions and proceed to identify the spline approximation that satisfies 
\[
K\mathbf{c} = f,
\]
with the matrices \( K, \mathbf{c}, \) and \( f \) explained in section 3.1. We define \( \{(\eta_i)\}_{i=1}^{n_b} \) as a collection of points situated on the boundary \( \partial \Omega \). For a given \( \lambda > 0 \), we seek \( \mathbf{c} \) that minimizes
\[
J(\mathbf{c}) = \sum_{i=1}^N (K_i\mathbf{c} - f_i)^2 + \lambda \sum_{i=1}^{n_b} [s_u(\eta_i) - g(\eta_i)]^2,
\]
where \( n_b \) denotes the number of boundary points and \( K_i, f_i \) represent the ith row of the matrix \( K, f \) respectively, and \( s_u = (s_1, \ldots, s_d) \). This minimization problem results in a system of linear equations.

\begin{example}
	To assess the performance of our method, we offer examples of solving both the 2D Stokes equations and Navier-Stokes equation within the domain illustrated in Figure \ref{curveddomain}, maintaining \( \lambda = 1 \) across all cases. We evaluate \( (\mathbf{u}^1, p^1) \) in Example \ref{ex1stokes} over the domains displayed in Figure \ref{curveddomain}. The Figures \ref{figure:cstokes1} to \ref{figure:cstokes4} showcase the spline approximations for the 2D Stokes equations in the designated domains. The associated numerical errors are tabulated in Table \ref{table:curvedStokes1}, underscoring the effectiveness of our method in handling curved domains, particularly when \( \mu= 1\text{e}-06 \).
	
	In addition to the 2D Stokes equations, we applied our method to the Navier-Stokes equation utilizing the same function. The outcomes, summarized in Table \ref{table:curvedStokes2}, exhibit comparable effectiveness, highlighting the versatility and consistency of our approach in the context of the Navier-Stokes equation with \( \mu=1 \).
	
	\begin{figure}[htbp]
		\centering
		\includegraphics[width=\linewidth]{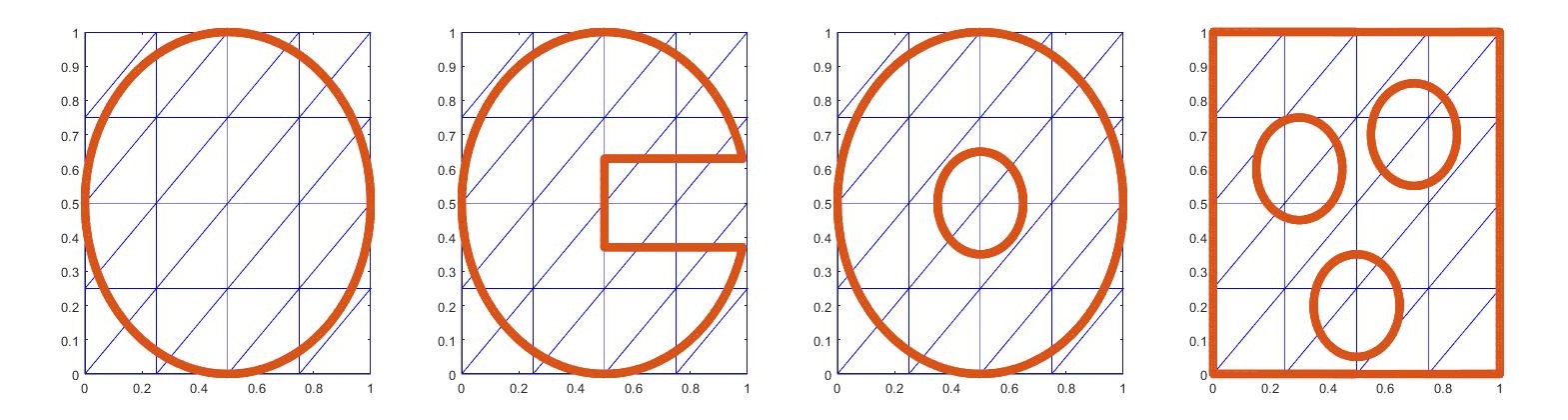}
		\caption{Domains defined by a curve immersed in mesh with $h=1/4$ }\label{curveddomain}
	\end{figure} 
	\begin{figure}[htbp]
		\centering
		\includegraphics[width=\linewidth]{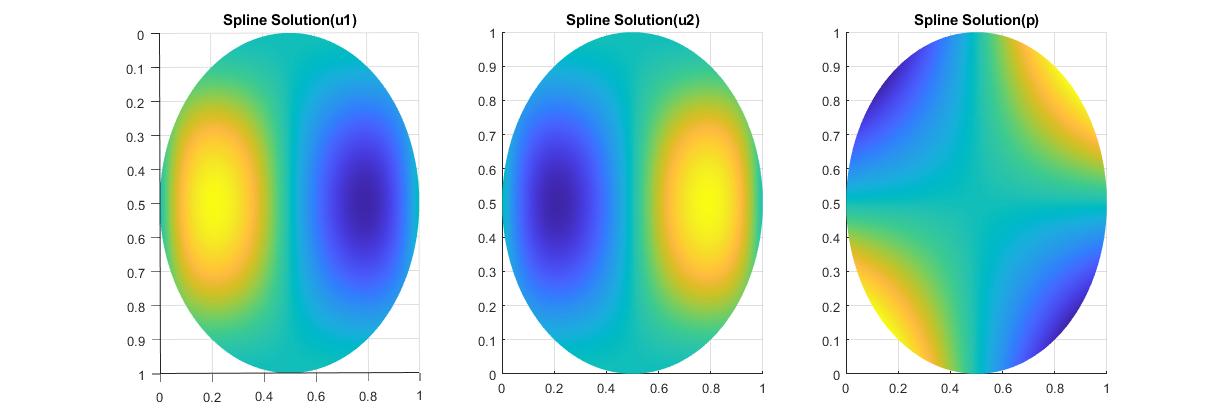}
		\caption{Case 1: Spline solutions $u_1, u_2, p$ with $D=7, r=2, \mu=1$ }\label{figure:cstokes1}
	\end{figure}
	\begin{figure}[htbp]
		\centering
		\includegraphics[width=\linewidth]{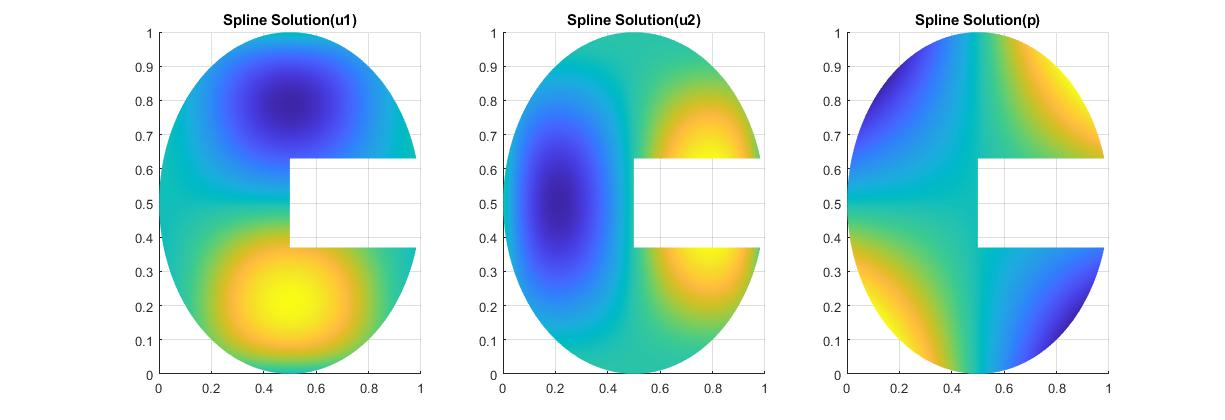}
		\caption{Case 2: Spline solutions $u_1, u_2, p$ with $D=7, r=2, \mu=1$ }
	\end{figure}
	\begin{figure}[htbp]
		\centering
		\includegraphics[width=\linewidth]{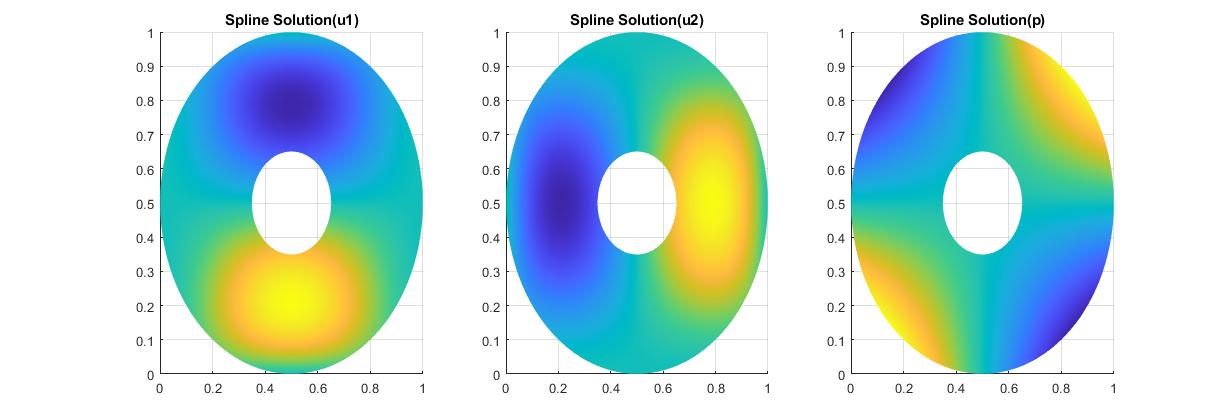}
		\caption{Case 3: Spline solutions $u_1, u_2, p$ with $D=7, r=2, \mu=1$ }
	\end{figure}
	\begin{figure}[htbp]
		\centering
		\includegraphics[width=\linewidth]{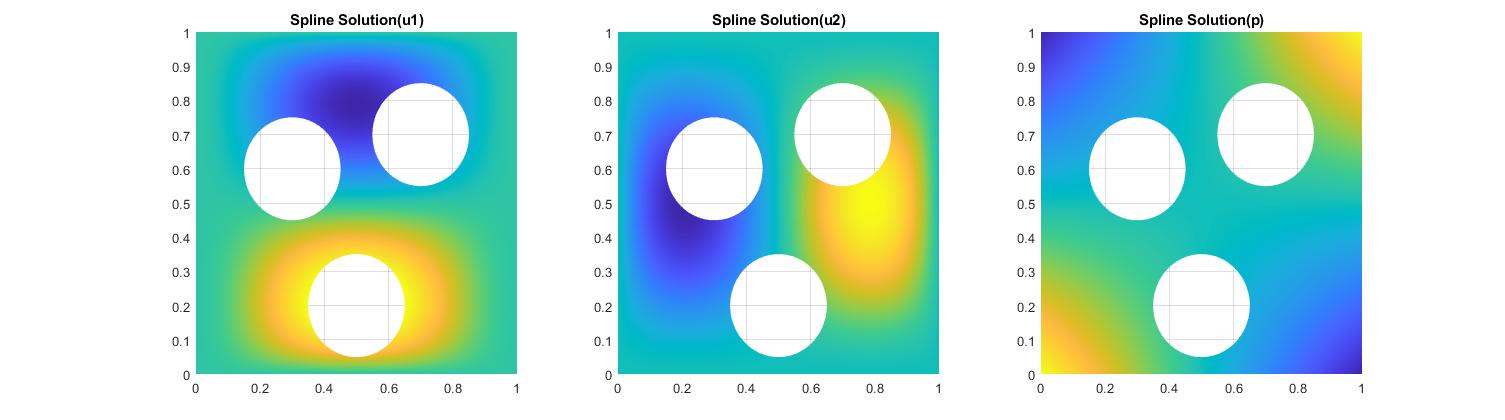}
		\caption{Case 4: Spline solutions $u_1, u_2, p$ with $D=7, r=2, \mu=1$ }\label{figure:cstokes4}
	\end{figure}
	\begin{table}[htp!]
		\centering
		\small
		\begin{tabular}{ c |c c c c c} 
			\hline
			case&CPU time&$\|\nabla \cdot\textbf{u}\|_\infty$&$|e_s|_{l_2} $&$|e_s|_{h_1} $& $|p-p_s|_{l_2} $\\
			\hline
			Case 1	&0.36  & 4.55e-14 & 3.37e-11 & 2.34e-10 & 1.17e-15   \\   
			Case 2&	0.38  & 2.11e-14 & 2.18e-12 & 2.70e-11 & 1.33e-15   \\   
			Case 3&	0.36  & 2.43e-14 & 2.82e-12 & 2.68e-11 & 1.58e-15   \\   
			Case 4&	0.36  & 3.95e-14 & 2.77e-13 & 5.69e-12 & 1.79e-15   \\   
			\hline
		\end{tabular}
		\caption{The $l_2$ errors and the divergence norm of $e_s, p-p_s$ for the spline solution of $(\textbf{u}^{s1},p^{s1})$ for the 2D Stokes equations over several curved boundary domains. Computations are carried out using the LL method with parameters $h=1/4,\mu=1e-06$, and $r=2$ for each degree $D=7$.}\label{table:curvedStokes1}
	\end{table}
	
	\begin{table}[htp!]
		\centering
		\small
		\begin{tabular}{ c |c c c c c} 
			\hline
			case&CPU time&$\|\nabla \cdot\textbf{u}\|_\infty$&$|e_s|_{l_2} $&$|e_s|_{h_1} $& $|p-p_s|_{l_2} $\\
			\hline
			Case 1& 0.89 & 3.78e-15 & 1.95e-15 & 2.39e-14 & 3.61e-14   \\        
			Case 2& 0.83 & 8.56e-14 & 1.78e-14 & 2.64e-13 & 2.15e-15   \\   
			Case 3& 0.94 & 7.72e-14 & 1.30e-14 & 2.22e-13 & 2.02e-15   \\   
			Case 4& 0.49 & 2.19e-15 & 1.09e-15 & 1.19e-14 & 2.49e-14 \\ 
			\hline
		\end{tabular}
		\caption{The $l_2$ errors and the divergence norm of $e_s, p-p_s$ for the spline solution of $(\textbf{u}^{s1},p^{s1})$ for the 2D Navier-Stokes equations over several curved boundary domains. Computations are carried out using the LL method with parameters $h=1/4,\mu=1$, and $r=2$ for each degree $D=7$.}\label{table:curvedStokes2}
	\end{table}
\end{example}
\begin{example}
	Next, we solve the Navier-Stokes equation on $[0,1]^2$ with $\mu=1$ and the exact functions are chosen as follows: 
	\[
	\begin{pmatrix}
		\textbf{u}^{c1}(x,y) \\
		p^{c1}(x,y) 
	\end{pmatrix} = 
	\begin{pmatrix}
		2\pi \sin^2(n\pi x) \cos(n\pi y) \sin(n\pi y) \\
		-2\pi \sin(n\pi x) \sin(n\pi y) \cos(n\pi x) \sin(n\pi y) \\
		\cos(n\pi x) \cos(n\pi y)
	\end{pmatrix}
	\]
	where $n=2.$ We refine the size $h$ of mesh which includes the curved domain $\Omega$. 
	The computed $l_2, h_1$ errors and convergence rates between the spline approximation and the exact solutions are tabulated in Table \ref{table:cNStokesconv}. This result shows that the our method work very well over the curved domain in Case 1.
	\begin{table}[htp!]
		\centering
		\small
		\begin{tabular}{ c | c c c c c c } 
			\hline
			$h$& $|e_s|_{l_2}$ & Rate & $|e_s|_{h_1}$ & Rate & $|p-p_s|_{l_2}$& Rate \\
			\hline
			$1/2$	&3.59e-01& -  & 6.30e+00 & - & 1.46e+01 & -   \\   
			$1/4$	&3.36e-03& 6.74  & 6.33e-02 & 6.64 & 2.38e-01 & 5.94   \\   
			$1/8$	&2.76e-04& 3.61  & 7.36e-03 & 3.10 & 1.61e-02 & 3.89   \\   
			$1/16$	&5.19e-06& 5.73  & 1.31e-04 & 5.81 & 2.71e-04 & 5.89   \\    
			\hline
		\end{tabular}
		\caption{$l_2, h_1$ errors and convergence rates for the spline solution of $(\textbf{u}^{c1},p^{c1})$ on the domain in Case 1 from the LL method with parameters $D=7$, and $r=2$ }
		\label{table:cNStokesconv}
	\end{table}
\end{example}

\section*{Acknowledgments}
We would like to acknowledge the assistance of volunteers in putting
together this example manuscript and supplement and anonymous referees for their valuable comments.
\bibliographystyle{siamplain}
\bibliography{references}

\end{document}